\let\ORIlabel\label
\let\ORIrefstepcounter\refstepcounter
\AddToHook{package/hyperref/before}{
   \let\label\ORIlabel
   \let\refstepcounter\ORIrefstepcounter}%

\documentclass[final,onefignum,onetabnum,letterpaper]{siamonline250211}

\usepackage[scale=0.8]{geometry} %
\usepackage[utf8]{inputenc} %
\usepackage[T1]{fontenc}    %

\usepackage{lipsum}
\usepackage{amsfonts}
\usepackage{thmtools}
\usepackage{epstopdf}
\usepackage{epstopdf-base}

\usepackage{datetime}
\usepackage{academicons}
\usepackage{xcolor}

\usepackage{amsmath}
\usepackage{amssymb}
\usepackage{commath}
\usepackage{mathtools}
\usepackage{bbm}

\usepackage{color}
\usepackage{graphicx}
\usepackage{grfext}

\usepackage{relsize}
\usepackage{adjustbox}
\usepackage{algorithm}
\usepackage[noend]{algpseudocode}
\usepackage{booktabs}
\usepackage{tikz}

\usepackage{verbatim}

\usepackage{amsopn}

\usepackage{todonotes}

\usepackage{lineno}
\usepackage{enumitem}

\usepackage{hyperref}
\usepackage[capitalize]{cleveref}

\usepackage{scalerel,stackengine}
\stackMath
\newcommand\what[1]{%
\ThisStyle{\savestack{\tmpbox}{\stretchto{%
  \scaleto{%
    \scalerel*[\widthof{\ensuremath{\SavedStyle #1}}]{\kern2pt\mathchar"0362\kern2pt}%
    {\rule{0ex}{\textheight}}%
  }{\textheight}%
}{1.2ex}}%
\stackon[-3pt]{\SavedStyle #1}{\tmpbox}%
}}

\ifpdf
    \DeclareGraphicsExtensions{.eps,.pdf,.png,.jpg}
\else
  \DeclareGraphicsExtensions{.eps}
\fi

\newdateformat{monthyeardate}{%
  \monthname[\THEMONTH] \THEDAY, \THEYEAR}

\renewcommand{\orcid}[1]{\href{https://orcid.org/#1}{\textcolor[HTML]{A6CE39}{orcid.org/#1}}}

\allowdisplaybreaks

\setlist[enumerate]{leftmargin=.5in}
\setlist[itemize]{leftmargin=.5in}

\newsiamthm{problem}{Problem}
\newsiamremark{remark}{Remark}
\newsiamremark{hypothesis}{Hypothesis}
\crefname{hypothesis}{Hypothesis}{Hypotheses}
\newsiamremark{example}{Example}
\newsiamthm{claim}{Claim}
\newsiamthm{conjecture}{Conjecture}

\headers{
Regularity-informed data assimilation
}{J.\ Glaubitz, D.\ Sharp, M.\ Le Provost, and Y.\ Marzouk}

\title{%
Regularity-informed data assimilation: A hierarchical Bayesian approach to ensemble Kalman filtering for hyperbolic conservation laws%
\thanks{
\monthyeardate\today
\corresponding{Jan Glaubitz}
}}

\author{
Jan Glaubitz\thanks{
Division of Applied Mathematics,
Link\"oping University,
Sweden
(\email{jan.glaubitz@liu.se},
\orcid{0000-0002-3434-5563})
}
\and
Daniel Sharp\thanks{
Center for Computational Science and Engineering,
Massachusetts Institute of Technology, USA
(%
\email{dannys4@mit.edu}, \orcid{0000-0002-0439-5084}; %
\orcid{0000-0003-0396-5740}; %
\email{ymarz@mit.edu}, \orcid{0000-0001-8242-3290}%
)
}
\and
Mathieu Le Provost\footnotemark[3]
\and
Youssef M.\ Marzouk\footnotemark[3]
}

\definecolor{english}{rgb}{0.0, 0.5, 0.0}
\definecolor{blue1}{RGB}{0,76,153}
\definecolor{green1}{RGB}{1,109,57}
\definecolor{red1}{RGB}{153,0,0}
\definecolor{orange1}{RGB}{255,127,80}

\DeclareMathOperator*{\argmin}{arg\,min}

\DeclareMathOperator{\diag}{diag}

\newcommand{\N}{\mathbb{N}}

\newcommand{\R}{\mathbb{R}}

\newcommand{\ns}{{n_{u}}}
\newcommand{\no}{{n_{b}}}
\newcommand{\nt}{{n_{s}}}
\newcommand{\nIAS}{n_{\operatorname{IAS}}}

\newif\ifrasterize
\rasterizefalse

\ifrasterize
\else
    \PrependGraphicsExtensions*{.pdf}
\fi

\makeatletter
\AddToHook{cmd/appendix/before}{\def\cref@section@alias{appendix}\def\cref@subsection@alias{appendix}}{\def\Cref@section@alias{Appendix}\def\Cref@subsection@alias{Appendix}}
\makeatother

\ifpdf
\hypersetup{
  pdftitle={Regularity-informed data assimilation},
  pdfauthor={Glaubitz et al}
}
\fi

\begin{document}

\maketitle
\begin{abstract}
We propose a novel regularity-informed filtering framework for data assimilation in the context of hyperbolic conservation laws and other time-dependent partial differential equations. We focus on systems whose states exhibit steep gradients and jump discontinuities.
While filtering is widely used to improve numerical simulations by incorporating observational data, traditional filtering methods lack awareness of the spatial regularity of states produced in these systems.
As a result, data assimilation often produces unphysical state estimates, introducing spurious oscillations in smooth regions and smearing sharp features.
To address this limitation, we introduce a filtering framework that incorporates edge-preserving regularization into the filter's analysis step; this framework balances simulation forecasts, observational data, and structural prior knowledge.
We formalize this approach using the ensemble Kalman filter (EnKF) and a class of hierarchical generalized sparse Bayesian learning (GSBL) priors, which adaptively infer spatially varying hyperparameters to promote non-oscillatory behavior in smooth regions while preserving discontinuities.
We demonstrate the effectiveness of the resulting GSBL-EnKF method on challenging benchmark problems governed by hyperbolic conservation laws.
Our results show that enforcing regularity in the analysis step yields sharper, less oscillatory state estimates and lower errors of the ensemble mean. This sometimes comes at the cost of ensemble spread, which we quantify and discuss.

\end{abstract}

\begin{keywords}
	Data assimilation,
    filtering,
    hyperbolic conservation laws,
    time-dependent PDEs,
    preserving regularity,
    generalized sparse Bayesian learning
\end{keywords}

\begin{AMS}
35L65,
35L67,
62F15,
62M20
\end{AMS}

\begin{Code}
    \url{https://github.com/dannys4/HierarchicalDAExamples/}
\end{Code}

\begin{DOI}
	Not yet assigned
\end{DOI}

\section{Introduction}
\label{sec:introduction}

Data assimilation (DA) is a potent framework for enhancing the accuracy of numerical simulations of complex physical systems by incorporating observational data and quantifying uncertainty~\cite{asch2016data,carrassi2018data,evensen2022data,sanz2023inverse}.
While there are many methods in DA, we consider the problem of \textit{ensemble filtering}, in which a set of identically distributed state estimates is simultaneously and continuously updated as new observational data become available over time.
Filtering techniques are widely applied in diverse scientific domains---e.g., weather prediction, systems biology, and robotics~\cite{evensen2022data,asch2016data}---to mitigate errors accumulated during simulation and to enhance predictive reliability.

While there are well-established filtering methods for diverse applications, faithfully reproducing physical phenomena in high-dimensional scientific settings is often challenging.
For these challenging applications, we consider the ensemble Kalman filter (EnKF) as the workhorse method for ensemble filtering.
One particularly demanding application is quantifying uncertainty for hyperbolic conservation laws, which describe a broad range of interesting dynamical systems~\cite{lax1973hyperbolic,dafermos2016hyperbolic,evans2022partial}.
Such systems have physical features that are often difficult to capture for filtering methods~\cite{bardos2005data,chapelle2012improving,boulanger2015data}.
This is because hyperbolic conservation laws lack inherent dissipation and their solutions can develop sharp gradients and shock discontinuities, even for smooth initial states and in finite time.
Traditional filtering methods, however, often produce ensembles lacking the piecewise regularity characteristic of hyperbolic systems.
Without delicate care, ensemble filters may therefore yield non-physical and unreliable results---often introducing spurious oscillations in smooth regions or smearing out discontinuities---unless carefully designed and constrained.

\subsection{Our contributions}

We propose a new framework for \textit{regularity-informed filtering}, tailored to hyperbolic conservation laws and, more broadly, to time-dependent partial differential equations (PDEs) whose solutions exhibit piecewise-smooth behavior.
That is, PDE solutions with a few jump discontinuities (or steep gradients) that behave smoothly between these irregularities.
While many modern numerical solvers are carefully designed to balance high accuracy in smooth regions with the ability to preserve discontinuous features~\cite{leveque1992numerical,leveque2002finite,shu2006essentially,hesthaven2007nodal,cockburn2012discontinuous,tadmor2012review,godlewski2013numerical,hesthaven2017numerical,glaubitz2020shock,offner2023approximation}, traditional ensemble filtering techniques often fail to do so.
Ultimately, the unphysical behavior found in many ensemble filtering methods degrades accuracy, obscures jump locations, and reduces overall reliability.

To address these issues, our proposed regularity-informed filtering framework introduces regularization into the filter's analysis step to promote the desired regularity characteristics.
We formalize our approach using the EnKF, whose analysis step can be interpreted as solving a randomized optimization problem with an objective inspired by Bayes rule.
Our central contribution is to \textit{tilt} the traditional EnKF prior in this Bayesian setup using an edge-preserving prior that encodes the belief that the true state is piecewise smooth.
We construct this prior using the edge-preserving generalized sparse Bayesian learning (GSBL) framework~\cite{calvetti2020sparse,glaubitz2023generalized,glaubitz2026efficient}, which employs a hierarchical approach to promote sparsity of discontinuities within and between samples.
Our GSBL approach thus encourages the state to be sharp at discontinuities (which are rare) and less oscillatory in smooth regions (which are common).
For efficient maximum a posteriori (MAP) estimation in the resulting \emph{GSBL-EnKF} posterior, we use a block coordinate descent approach, which is effective for linear inverse problems with GSBL priors~\cite{calvetti2020sparse,calvetti2020sparsity,xiao2023sequential,glaubitz2023leveraging,lindbloom2024generalized,lindbloom2025priorconditioned}.
Notably, our approach is developed and used for state discretizations on arbitrary (potentially non-uniform) Cartesian grids, in contrast to many existing approaches that typically assume equidistant grids.
The sparsifying transform is built from the differentiation and quadrature operators of the numerical discretization itself, so that the approach applies directly to the non-equidistant nodal sets of high-order discontinuous Galerkin methods.
In this work, we consider one-dimensional discretizations and discuss the extension to multiple dimensions in \Cref{sub:problems}.

We demonstrate the performance of the proposed GSBL-EnKF approach on several challenging test problems, including the one-dimensional linear advection, the Burgers equation, and the compressible Euler equations.
We find that the GSBL approach can substantially improve upon the typical EnKF, especially as the observation grid becomes denser, and our method can quickly converge to a reasonable ensemble even with adversarial initializations.
Our results indicate that preserving regularity during the analysis step yields ensemble members that are piecewise smooth in space and an ensemble mean with lower error.
At the same time, the regularized update reduces the ensemble spread, and we quantify this effect on the uncertainty quantification (UQ) provided by the filter.
While the present study focuses on hyperbolic conservation laws, the EnKF, and GSBL priors, the broader principle---incorporating problem-specific regularity assumptions into the DA process---is more widely applicable.
We therefore expect this work to pave the way for further regularity-informed DA techniques across a range of dynamical systems and filtering paradigms.

\subsection{Related work}

Despite its importance in many applications, relatively few studies have focused on preserving regularity in state estimation problems, especially for piecewise-smooth profiles with jump discontinuities.
Some early contributions~\cite{freitag2010l1,budd2011regularization,freitag2013resolution} revisit the four-dimensional variational approach (4DVar) from an inverse-problem perspective, showing that 4DVar can be interpreted as a form of Tikhonov regularization.
Drawing inspiration from image processing, the studies introduce both $L_1$ and mixed $L_1$--$L_2$ regularizations into 4DVar.
A fundamentally different approach is proposed in~\cite{srivastava2023feature}, introducing feature-informed DA.
In this method, feature information---such as shock locations or wavefronts---is incorporated via set-valued observation operators.
These operators involve a search over the state space, making the approach qualitatively distinct from traditional observation models.
In~\cite{li2024structurally,li2025structurally,li2025structurally2}, structurally-informed priors are constructed for the ensemble transform Kalman filter method using second-moment information of state gradients.
The approach modifies the prior weighting matrix in the assimilation objective by clustering gradient statistics to emphasize regions of discontinuity, thereby replacing standard covariance models.

Closest to the present work is~\cite{kim2023hierarchical}, which incorporates sparsity-promoting regularization into ensemble Kalman methods for static inverse problems $\mathbf{b} = \mathcal{G}(\mathbf{u}) + \boldsymbol{\varepsilon}$: a conditionally Gaussian prior on $\mathbf{u}$ with generalized gamma hyper-priors is combined with an iterative ensemble Kalman method for the update of $\mathbf{u}$, alternating with closed-form updates of the hyperparameters to approximate the MAP estimate.
The present work differs in three respects.
First, we consider time-dependent PDEs and sequential filtering, where the forecast ensemble enters each analysis step as a Gaussian factor centered at the individual members, whereas~\cite{kim2023hierarchical} treats static problems and uses the ensemble as a derivative-free optimizer initialized from the prior; regularized ensemble Kalman filtering for DA is listed there as future work.
Second, in~\cite{kim2023hierarchical} a single hyperparameter vector is shared by the whole ensemble and updated from the ensemble mean, whereas here each ensemble member carries its own hyperparameter vector, so that the ensemble also represents uncertainty about the location of discontinuities.
Third, sparsity is imposed in~\cite{kim2023hierarchical} on the unknown itself, whereas we impose it on a transformed state $S\mathbf{u}$ with a discretization-dependent transform $S$, for which we derive a Kalman-type update with an augmented observation operator (\Cref{lem:Kalman}).
Regularization of ensemble Kalman methods has also been studied in the inverse-problem setting, including hierarchical parameterizations~\cite{chada2018parameterizations}, $\ell_p$-type sparsity regularization~\cite{lee2021lp,schneider2022ensemble}, and adaptive Tikhonov regularization~\cite{iglesias2021adaptive}; sparsity-promoting regularization of variational DA in transform domains is considered in~\cite{ebtehaj2014variational}.

A comparative study of sequential filtering methods in compressible flows with shock discontinuities is presented in~\cite{houba2024sequential}.
This work shows that both the EnKF and the extended Kalman filter may produce non-physical trajectories and values, such as negative pressures or irregular states, which can cause the numerical solver to break down.
The study in~\cite{edoh2025sequential} addresses the issue of non-physical values in DA by operating on logarithmically transformed variables for forecast updates, an example of Gaussian anamorphosis~\cite{amezcuaGaussianAnamorphosisAnalysis2014}.
To address non-Gaussianity, the authors of~\cite{hansen2024normal} propose a normal-score EnKF variant and apply it to systems with shocks.
A variational formulation tailored to shocks is proposed in~\cite{west2025variational}, where 4DVar is applied to 1D compressible flow with discontinuities.
The study carefully examines the feasibility and performance of adjoint-based methods for shock-resolving DA.
The work~\cite{subrahmanya2025feature} develops a feature-preserving ensemble transform particle filter for the compressible Euler equations.
By combining particle-based DA with sequence alignment techniques, the method aligns and interpolates particles along feature-consistent paths.

Lastly, recent work adapts the EnKF to shock-laden flows using machine learning.
In~\cite{zhou2026neural}, the authors argue that the EnKF performs poorly for discontinuous states due to a bimodal forecast distribution that can arise near jump discontinuities; this bimodality violates the Gaussian assumption underpinning the EnKF.
The authors map the forecast ensemble to the parameter space of a neural network (i.e., weights and biases), where they perform DA.
Notably, this approach requires training a neural network for each ensemble member, at each analysis step, and for each PDE component, which they do in sequence using a nearest-neighbor ordering of members (see~\cite[Figure 5]{zhou2026neural}).
A similar approach is taken in~\cite{chandravamsi2026feature}, which performs a feature-preserving ensemble update in a learned low-dimensional latent space, where shock and flow features admit a smooth manifold
representation.
Using an encoder-decoder architecture eliminates the need for sequential training, enabling the authors to map the ensemble to physical space via a shared decoder.

The above contributions collectively highlight the challenges in constructing DA algorithms that preserve sharp features.
They also emphasize the need for further methodological innovations that balance computational efficiency with physical realism in nonlinear systems.

\subsection*{Outline}

\Cref{sec:prelim} reviews essential background on DA, filtering, and the EnKF.
\Cref{sec:proposed} introduces the proposed regularity-informed filtering framework.
\Cref{sec:tests} evaluates the performance of the resulting GSBL-EnKF method on several challenging benchmark problems governed by hyperbolic conservation laws.
\Cref{sec:summary} provides a summary, articulates open problems, and discusses future research directions.

\section{Preliminaries}
\label{sec:prelim}

We provide some preliminaries on the general problem setup, filtering, the EnKF, and hyperbolic conservation laws.
This section is mainly based on the monographs~\cite{asch2016data,carrassi2018data,evensen2022data,sanz2023inverse}.

\subsection{Problem setup and notation}
\label{sub:prelim_problem}

Consider a domain $\Omega \subset \R^d$ and a state variable $u(x,t)$, which depends on time $t \in [0,T]$ and spatial position $x \in \Omega$.
Furthermore, let $\mathbf{u}_j \in \R^\ns$, $j=0,\dots,J$, be discrete numerical approximations of $u(x,t_j)$ at some grid points $\{ x_i \}_{i=1}^{\ns} \subset \Omega$, i.e., $(\mathbf{u}_j)_i \approx u(x_i,t_j)$.
We assume that the discrete states $\mathbf{u}_j$ follow a given dynamics model:
\begin{equation}\label{eq:dynamics_model}
	\mathbf{u}_{j} = \Psi( \mathbf{u}_{j-1} ) + \boldsymbol{\xi}_{j}, \quad j=1,\dots,J,
\end{equation}
where $\Psi: \R^\ns \to \R^\ns$ is a known, time-independent, solution operator and $\boldsymbol{\xi}_j$ is independent and identically distributed (i.i.d.) noise with $\boldsymbol{\xi}_j \sim \mathcal{N}( \mathbf{0}, \Gamma_x )$.
In practice, the operator $\Psi$ often corresponds to simulating a time-dependent PDE using a numerical solver.
In this case, the additive noise $\boldsymbol{\xi}_{j}$ can model the approximation errors introduced by the numerical solver in addition to inherent stochasticity in the system.

We use the above assumptions to model our incomplete knowledge of the system's true initial state and dynamics, which is a typical challenge in many applications.
Errors and uncertainties of the initial state $\mathbf{u}_0$ consequently propagate through the model~\cref{eq:dynamics_model} into the estimates for the later states $\mathbf{u}_j$, making reliable predictions difficult.
Fortunately, in many cases, we have access to noisy observations of the true state, which can be used to correct errors in the estimates $\mathbf{u}_j$.
We formalize this via the observation model:
\begin{equation}\label{eq:data_model}
	\mathbf{b}_{j} = H \mathbf{u}_{j}  + \boldsymbol{\eta}_{j}, \quad j=1,\dots,J,
\end{equation}
where $\mathbf{b}_j \in \R^\no$ is the observational data, $H \in \R^{\no \times \ns}$ is a known linear observation operator, and $\boldsymbol{\eta}_j$ is additive observation noise.
We make the common assumption that the observational noise sequence $\{\boldsymbol{\eta}_j\}$ is i.i.d.\ with $\boldsymbol{\eta}_j \sim \mathcal{N}( \mathbf{0}, \Gamma )$ and independent of $\mathbf{u}_0$. As a point of notation, we write $\mathbf{z}_{i:j}$ for the time series $\mathbf{z}_{i},\mathbf{z}_{i+1},\dots,\mathbf{z}_{j}$.

\subsection{Filtering}
\label{sub:prelim_filtering}

Filtering refers to sequentially updating the probability distribution of the state $\mathbf{u}_j$, denoted $\pi( \mathbf{u}_j | \mathbf{b}_{1:j} )$, when data $\mathbf{b}_j$ is acquired.
At time $t_j$, the \textit{filtering distribution} $\pi( \mathbf{u}_{j} | \mathbf{b}_{1:j} )$ is typically derived from the previous time's distribution $\pi( \mathbf{u}_{j-1} | \mathbf{b}_{1:j-1} )$ through a two-step process:
First, in the \emph{forecast step}, we create samples distributed according to $\pi( \mathbf{u}_{j} | \mathbf{b}_{1:j-1} )$ via simulating the dynamics model~\cref{eq:dynamics_model} for samples of the previous filtering distribution $\pi( \mathbf{u}_{j-1} | \mathbf{b}_{1:j-1} )$.
This is followed by the \emph{analysis step}, in which the observational data from~\cref{eq:data_model} is incorporated to create samples distributed (approximately) according to $\pi( \mathbf{u}_{j} | \mathbf{b}_{1:j} )$.
To this end, the density of the posterior distribution $\pi( \mathbf{u}_{j} | \mathbf{b}_{1:j} )$ can be obtained via Bayes' rule:
\begin{equation}\label{eq:filtering_Bayes}
	\pi( \mathbf{u}_{j} | \mathbf{b}_{1:j} )
		\propto \pi( \mathbf{b}_{j} | \mathbf{u}_{j}, \mathbf{b}_{1:j-1} ) \,
		\pi( \mathbf{u}_{j} | \mathbf{b}_{1:j-1} ) = \pi( \mathbf{b}_{j} | \mathbf{u}_{j} ) \,
		\pi( \mathbf{u}_{j} | \mathbf{b}_{1:j-1} ).
\end{equation}
We encode our assumptions about the data model~\cref{eq:data_model} with $\pi( \mathbf{b}_{j} | \mathbf{u}_{j} )$,%
\footnote{Recall that $\mathbf{b}_{j}$ is independent of $\mathbf{b}_{1:j-1}$ given $\mathbf{u}_j$, seen in~\cref{eq:data_model}.}
while the distribution $\pi( \mathbf{u}_{j} | \mathbf{b}_{1:j-1} )$ serves as a prior, incorporating our a priori assumptions about the state $\mathbf{u}_{j}$.
In the case of filtering, such a priori assumptions include that $\mathbf{u}_{j}$ is generated from $\mathbf{u}_{j-1}$ in the prediction step through the dynamics model~\cref{eq:dynamics_model}.

\subsection{Ensemble Kalman filtering}
\label{sub:prelim_EnKF}

There are several filtering methods for performing the above prediction and analysis steps.
Here, we focus on the EnKF, which uses an ensemble of $P$ particles $\{ \mathbf{u}_j^{(p)} \}_{p=1}^P$ to inform the filtering update.
The particles are all given equal weight, making it possible to approximate the filtering distribution via $\pi( \mathbf{u}_{j} | \mathbf{b}_{1:j} ) \approx \frac{1}{P} \sum_{p=1}^P \delta( \mathbf{u}_{j} - \mathbf{u}_j^{(p)} )$.
This approximation can be accurate if $ P$ is sufficiently large and the filtering distributions are approximately Gaussian; see~\cite[Section 2.7]{calvelloEnsembleKalmanMethods2025} and references therein.
In problems where approximate Gaussianity of the filtering distribution fails---for instance, due to strong nonlinearity of $\Psi$ and large observation noise---the EnKF is perhaps better understood as a sequential optimization method, similar in spirit to 3DVar~\cite{lorenc1986analysis,abarbanel2013predicting,sanz2023inverse,calvelloEnsembleKalmanMethods2025}.

Suppose we are given an ensemble of initial states $\mathbf{u}_{0}\sim\pi_0$ for some distribution $\pi_0$.
To generate the ensemble at time $j$ from the previous one, we proceed as follows.
First, in the forecast step, we propagate each particle through the dynamics model~\cref{eq:dynamics_model} to obtain
\begin{equation}\label{eq:EnKF_prediction}
	\what{\mathbf{u}}_{j}^{(p)} \sim\mathcal{N}\left( \Psi\left( \mathbf{u}_{j-1}^{(p)} \right), \Gamma_x \right), \quad p=1,\dots,P.
\end{equation}
We then compute the sample prior mean and covariance matrix as
\begin{equation}\label{eq:EnKF_prior_covariance}
	\what{\mathbf{m}}_{j}
		= \frac{1}{P} \sum_{p=1}^P \what{\mathbf{u}}_{j}^{(p)}, \quad
	\what{C}_{j}
		= \frac{1}{P-1} \sum_{p=1}^P \left( \what{\mathbf{u}}_{j}^{(p)} - \what{\mathbf{m}}_{j} \right) \left( \what{\mathbf{u}}_{j}^{(p)} - \what{\mathbf{m}}_{j} \right)^\top.
\end{equation}
Next, in the analysis step, we update the ensemble based on the observational data $\mathbf{b}_j$.
To this end, consider the set of simulated observations $\mathbf{b}_{j}^{(p)} = \mathbf{b}_{j} +  \boldsymbol{\eta}_{j}^{(p)}$ with i.i.d. noise $\boldsymbol{\eta}_{j}^{(p)} \sim \mathcal{N}(\mathbf{0}, \Gamma)$ for $p=1,\ldots,P$.
We refer to the particles $\mathbf{b}_{j}^{(p)}$ as \textit{perturbed observations}.
Each particle is then updated according to the \emph{Kalman update}:
\begin{equation}\label{eq:EnKF_analysis1}
	\mathbf{u}_{j}^{(p)}
		= \what{\mathbf{u}}_{j}^{(p)}
			+ K_{j} \left( \mathbf{b}_{j}^{(p)} - H \what{\mathbf{u}}_{j}^{(p)} \right), \quad
    K_{j} = \what{C}_{j} H^\top \left( H \what{C}_{j} H^\top + \Gamma \right)^{-1},
\end{equation}
where $K_{j}$ is the so-called \emph{Kalman gain}.
Notably, the analysis step~\cref{eq:EnKF_analysis1} can be re-written as a quadratic minimization problem:
\begin{equation}\label{eq:EnKF_analysis2}
	\mathbf{u}_{j}^{(p)}
		= \argmin_{\mathbf{u} \in \R^\ns} \left\{
			\norm{ H \mathbf{u} - \mathbf{b}_{j}^{(p)} }^2_{\Gamma} +
			\norm{ \mathbf{u} - \what{\mathbf{u}}_{j}^{(p)} }^2_{\what{C}_{j}}
		\right\},
\end{equation}
where the above norms are defined through $\| \mathbf{v} \|^2_{A} = \| A^{-1/2} \mathbf{v} \|_2^2$ for a symmetric positive definite $A$.
For more information, see~\cite{sanz2023inverse} and references therein.
In \Cref{sec:proposed}, we use~\cref{eq:EnKF_analysis2} as a starting point for developing the proposed regularity-informed DA framework. In particular, each individual ensemble member produced by~\cref{eq:EnKF_analysis2} can be interpreted as the MAP estimator of the posterior:
\begin{equation}\label{eq:randomized_MAP}
	\pi(\mathbf{u}^{(p)}_j|\mathbf{b}^{(p)}_j,\what{\mathbf{u}}_j^{(p)}) \propto \pi(\mathbf{b}^{(p)}_j|\mathbf{u}^{(p)}_j)\pi(\mathbf{u}^{(p)}_j|\what{\mathbf{u}}_j^{(p)}).
\end{equation}
Here, the forecast member $\what{\mathbf{u}}_j^{(p)}$ plays the role of the prior mean and the sample covariance $\what{C}_j$ that of the prior covariance, i.e., $\pi(\mathbf{u}^{(p)}_j|\what{\mathbf{u}}_j^{(p)}) = \mathcal{N}(\what{\mathbf{u}}_j^{(p)}, \what{C}_j)$. This is the randomized maximum likelihood (RML) view of the perturbed-observation EnKF~\cite{kitanidis1995quasi,oliver1996conditioning,chen2012ensemble,bardsley2014randomize}: For a linear observation model and Gaussian prior, the minimizers of the randomized objectives~\cref{eq:EnKF_analysis2} are exact samples from the Gaussian posterior, whereas for the regularized (and hence nonlinear) updates considered below they are only approximate samples.
In EnKF implementations, practitioners often adopt various regularization heuristics for $\what{C}_j$, e.g., localization and inflation.
These will be discussed further in \Cref{sub:problems}, and our choices are described in \Cref{app:details}.

\subsection{Hyperbolic conservation laws}
\label{sub:prelim_CLs}

We focus on dynamical models~\cref{eq:dynamics_model} corresponding to time-dependent PDEs whose solutions are piecewise smooth and may exhibit a finite number of jump discontinuities.
A particularly important subclass of such problems is hyperbolic conservation laws, which arise in modeling a wide range of physical phenomena~\cite{lax1973hyperbolic,dafermos2016hyperbolic,evans2022partial}.
A generic system of conservation laws can be written as
\begin{equation}
	\partial_t \boldsymbol{u}(x,t) + \nabla \cdot \boldsymbol{f}( \boldsymbol{u}(x,t) ) = 0, \quad x \in \Omega \subset \R^d, \ t > 0,
\end{equation}
where $\boldsymbol{u}: \Omega \times \R_{>0} \to \R^m$ is a vector-valued function with $m$ so-called conserved variables and $\boldsymbol{f}: \R^m \to \R^{m \times d}$ is a smooth flux function, representing the flow of these conserved variables;
in this work, we focus on dimension $d = 1$.
Notable examples of conservation laws include the Burgers, shallow water, and Euler equations.

Solutions of hyperbolic conservation laws can develop jump discontinuities, even in finite time and for smooth initial states---a well-known observation~\cite{riemann1860fortpflanzung} that leads to significant theoretical and practical challenges.
Numerical methods for hyperbolic conservation laws typically leverage physical knowledge (e.g., conservation, entropy, and upwinding considerations) to address these challenges.
Further, many methods are designed to capture discontinuities without introducing spurious oscillations, i.e., to keep the numerical solution piecewise smooth.
Exploiting this type of regularity has led to the design of highly efficient solvers for conservation laws;
see~\cite{leveque1992numerical,leveque2002finite,shu2006essentially,hesthaven2007nodal,cockburn2012discontinuous,tadmor2012review,godlewski2013numerical,hesthaven2017numerical,glaubitz2020shock,offner2023approximation} and references therein.

\section{Proposed framework: regularity-informed filtering}
\label{sec:proposed}

We propose \emph{GSBL-EnKF}---a new DA framework modifying the EnKF's analysis step, motivated by piecewise-smooth solutions to hyperbolic conservation laws.
Henceforth, we consider the analysis step at a generic time $t_{j}$.
For simplicity, we thus drop the subindex ``$j$'' and write $\mathbf{u}^{(p)}$, $\mathbf{b}^{(p)}$, $\what{\mathbf{u}}^{(p)}$, and $\what{C}$ instead of $\mathbf{u}^{(p)}_{j}$, $\mathbf{b}^{(p)}_{j}$, $\what{\mathbf{u}}^{(p)}_{j}$, and $\what{C}_{j}$, respectively.

\subsection{The basic idea}\label{sec:basic_idea}

Our idea for promoting piecewise-smooth solutions in the EnKF's analysis step can be formalized by augmenting~\cref{eq:EnKF_analysis2} with an additional edge-preserving regularization term, yielding:\footnote{We use $\theta^{-1}$ for the regularization parameter to draw a clearer connection to GSBL below.}
\begin{equation}\label{eq:regEnKF_analysis}
	\mathbf{u}^{(p)}
		= \argmin_{\mathbf{u} \in \R^\ns} \left\{
			\norm{ H \mathbf{u} - \mathbf{b}^{(p)} }^2_{\Gamma}
			+ \frac{1}{\lambda}\norm{ \mathbf{u} - \what{\mathbf{u}}^{(p)} }^2_{\what{C}}
			+ \frac{1}{\theta} \mathcal{R}( \mathbf{u} )
		\right\}, \quad p=1,\dots,P.
\end{equation}
We interpret $\| H \mathbf{u} - \mathbf{b}^{(p)} \|^2_{\Gamma}$ as a \emph{data fidelity} term and $\| \mathbf{u} - \what{\mathbf{u}}^{(p)} \|^2_{\what{C}}$ as a \emph{forecast regularization} term.
The data fidelity term promotes the reconstruction $\mathbf{u}$ to reasonably explain the observed data, whereas the forecast regularization term ensures that the predictions align with the dynamics model~\cref{eq:dynamics_model}.
The motivation for $\mathcal{R}( \mathbf{u} )$ is that it acts as an additional, edge-preserving regularization term that promotes $\mathbf{u}$ to have a piecewise smooth profile.
Finally, the regularization parameters $\lambda,\theta > 0$ balance the trio of terms: data fidelity, forecast regularization, and edge-preserving regularization.

A natural---although potentially na\"ive---choice for $\mathcal{R}( \mathbf{u} )$ to promote piecewise-smooth solutions is the $\ell_1$-regularizer $\mathcal{R}( \mathbf{u} ) = \| S \mathbf{u} \|_1 = \sum_{k} | [S \mathbf{u}]_k |$,
where we choose an operator $S \in \R^{\nt \times \ns}$ and denote $[S\mathbf{u}]_k$ as the $k$-th entry of $S\mathbf{u} \in \R^\nt$.
This choice results in the $\ell^1$-regularized EnKF analysis step:
\begin{equation}\label{eq:l1EnKF_analysis}
	\mathbf{u}^{(p)}
		= \argmin_{\mathbf{u} \in \R^\ns} \left\{
			\norm{ H \mathbf{u} - \mathbf{b}^{(p)} }^2_{\Gamma}
			+ \frac{1}{\lambda}\norm{ \mathbf{u} - \what{\mathbf{u}}^{(p)} }^2_{\what{C}}
			+ \frac{1}{\theta} \| S \mathbf{u} \|_1
		\right\}, \quad p=1,\dots,P.
\end{equation}
In essence, the $\ell^1$-norm $\| \cdot \|_1$ acts as a convex surrogate for counting the number of nonzero entries, promoting $S \mathbf{u}$ to be sparse.
That is, we promote the vector $S\mathbf{u}$ to have only a few dominant entries, with the remaining being close to zero.
For example, if $S$ approximates a gradient operator, then $\mathcal{R}( \mathbf{u} ) = \| S \mathbf{u} \|_1$ promotes reconstructions that average out to few regions of large gradients, as these are sparse in the edge domain. For this reason, the user-specified map $S$ is often called a \textit{sparsifying transform}.
See \Cref{rem:sparsifying_transform} for details.

\begin{remark}
	The idea of incorporating $\ell_1$-regularization into the analysis step is inspired by~\cite{glaubitz2019high}, which considered~\cref{eq:l1EnKF_analysis} without the forecast term $\norm{ \mathbf{u} - \what{\mathbf{u}}^{(p)} }^2_{\what{C}}$.
	The resulting procedure served as a shock-capturing mechanism, used after each time step for a high-order discontinuous Galerkin method applied to hyperbolic conservation laws.
	In line with our motivation, the purpose of this work is to use the regularization term to remove spurious numerical oscillations around jump discontinuities.
\end{remark}

A common difficulty for regularized inverse problems such as~\eqref{eq:l1EnKF_analysis} is the selection of an appropriate regularization parameter $\theta$.
This parameter critically influences the quality of the reconstruction~\cite{vogel2002computational,hansen2010discrete}.
Furthermore, it is often argued that spatially varying regularization, $\mathcal{R}( \mathbf{u} ) = \sum_{k} \frac{1}{\theta_k} \left| [S \mathbf{u}]_k \right|$, should be used instead of a scalar parameter $\theta$ for promoting piecewise smooth reconstructions~\cite{candes2008enhancing,mansour2017recovery,adcock2019joint}.
To overcome the challenges of hyperparameter selection, we leverage a hierarchical Bayesian approach to inverse problems~\cite{kaipio2006statistical,stuart2010inverse,calvetti2023bayesian} to develop the proposed GSBL-EnKF methodology.
Two principles guide this construction: Structural knowledge about the state---here, piecewise smoothness---should be encoded in the prior.
At the same time, the resulting analysis step should remain about as tractable as the EnKF update it replaces. The $\ell_1$-type regularization may not achieve such tractability; it cannot generally be minimized analytically.
The hierarchical GSBL prior introduced below realizes the first principle, while the second is met because the associated MAP estimate can be computed using Kalman-type updates (cf. \Cref{sub:analysis_step}).

\begin{remark}\label{rem:sparsifying_transform}
	A typical choice for the sparsifying transform $S \in \R^{\nt \times \ns}$ is that of discrete gradient operators~\cite{rudin1992nonlinear,chan2000high,bredies2010total}, which frequently assume the nodal solution values $\mathbf{u}^{(p)}$ are given on equidistant grid points.
	Many modern high-order numerical PDE solvers, however, use non-equidistant grid points, or the vector $\mathbf{u}$ contains coefficients for a basis instead of interpolation values for the underlying state~\cite{kopriva2009implementing,chen2020review,ranocha2023efficient}.
    For this reason, we use $S$ as a scaled approximation of the second-derivative operator for the discretization. In particular, we use approximation $\|S\mathbf{u}_j\|_2^2\approx \int\left(\partial_x^2 u(x,t_j)\right)^2\,\mathrm{d}x$ for $[\mathbf{u}_j]_k \approx u(x_k,t_j)$.
    Here, the points $x_k$ are assumed to be the quadrature points on a discretized one-dimensional spatial grid and thus $S = \mathrm{diag}(\mathbf{w})^{1/2}D^2$, where the vector $\mathbf{w}\in\mathbb{R}^{\ns}$ collects positive quadrature weights, the map $D\in\mathbb{R}^{\ns\times\ns}$ is the correct differentiation operator for the numerical discretization, and the operator $\mathrm{diag}:\mathbb{R}^n\to\mathbb{R}^{n\times n}$ is, in this context, the mapping of a vector to a diagonal matrix.

    For the nodal discontinuous Galerkin discretization used in this paper, the operator $D$ is the element-wise differentiation matrix of the piecewise-polynomial approximation: the state values at the quadrature nodes of an element are mapped to the interpolation nodes, the polynomial is differentiated, and the result is evaluated at the quadrature nodes again.
    Consequently, the sparsification transform $S$ is block diagonal with one block per element, $\nt = \ns$, and its application $S\mathbf{u}$ measures the curvature of the state $\mathbf{u}$ within each element; jumps of the state across element interfaces are not penalized by the map $S$.
    In particular, states that are linear within each element lie in the kernel of $S$.
    We penalize the second rather than the first derivative so that profiles with a nonzero but constant slope are not penalized.
    A first-derivative transform would favor piecewise-constant reconstructions (staircasing), whereas the second-derivative transform favors profiles that are linear between the few nodes at which $[S\mathbf{u}]_k$ is large.
    The minimizer of~\cref{eq:l1EnKF_analysis} or any variant presented in our work will not quite achieve this piecewise, but nevertheless is penalized for deviating substantially, while the data-fidelity term of the problem mitigates dissipation.
    Finally, we note that the proposed GSBL-EnKF approach below is formally applicable to any sparsifying transform $S \in \R^{\nt \times \ns}$.
    For instance, we also experimented with polynomial annihilation operators~\cite{archibald2005polynomial,glaubitz2019high} for $S$, but found that using second derivatives improved on such maps.
    We leave a more thorough exploration of alternative choices of $S$ for promoting piecewise-smooth solutions to future work.
\end{remark}

\subsection{The GSBL-EnKF approach}

We now take a hierarchical Bayesian approach~\cite{calvetti2020sparse,kim2023hierarchical} and let hyperparameters $\boldsymbol{\theta}$ be uncertain and vary across the ensemble. Recalling the variational formulation of the analysis step~\cref{eq:EnKF_analysis2} and the randomized MAP estimation suggested by~\cref{eq:randomized_MAP}, we define the posterior of the $p$-th particle as:
\begin{equation}\label{eq:Bayes_rule}
	\pi( \mathbf{u}^{(p)}, \boldsymbol{\theta}^{(p)} | \mathbf{b}^{(p)},\what{\mathbf{u}}^{(p)})
		\propto \pi( \mathbf{b}^{(p)} | \mathbf{u}^{(p)}) \, \pi( \mathbf{u}^{(p)}, \boldsymbol{\theta}^{(p)} | \what{\mathbf{u}}^{(p)}),
\end{equation}
which combines the likelihood $\pi( \mathbf{b}^{(p)} | \mathbf{u}^{(p)} )$ implied by the data model~\cref{eq:data_model} with a prior $\pi( \mathbf{u}^{(p)}, \boldsymbol{\theta}^{(p)} | \what{\mathbf{u}}^{(p)})$ that encodes our structural beliefs about $\mathbf{u}^{(p)}$ given forecast $\what{\mathbf{u}}^{(p)}$.
Here, the set $\{\boldsymbol{\theta}^{(p)}\}_{p=1}^P$ is an ensemble of auxiliary hyperparameter vectors that remain to be specified below.
Denoting the collection of all particles by $\mathbf{u}^{(1:P)}$, the forecast ensemble by $\what{\mathbf{u}}^{(1:P)}$, and the perturbed observations by $\mathbf{b}^{(1:P)}$, the joint posterior density for the whole ensemble is given by:
\begin{equation}
	\pi( \mathbf{u}^{(1:P)}, \boldsymbol{\theta}^{(1:P)} | \mathbf{b}^{(1:P)},\what{\mathbf{u}}^{(1:P)}) = \prod_{p=1}^P \pi( \mathbf{u}^{(p)}, \boldsymbol{\theta}^{(p)} | \mathbf{b}^{(p)},\what{\mathbf{u}}^{(p)}).
\end{equation}
Below, we address the different components of the posterior~\cref{eq:Bayes_rule}.
Ultimately, the proposed GSBL-EnKF will correspond to the MAP estimate of the resulting posterior density.

\subsubsection{The likelihood}
\label{subsub:likelihood}

The likelihood $\pi( \mathbf{b}^{(p)} | \mathbf{u}^{(p)} )$ in~\cref{eq:Bayes_rule} models our assumptions about the data generating process~\cref{eq:data_model} and the noise characteristics.
We have the likelihood density $\pi(\mathbf{b}^{(p)} | \mathbf{u}^{(p)}) = \mathcal{N}( H \mathbf{u}^{(p)}, \Gamma )$, i.e.,
\begin{equation}\label{eq:likelihhod}
	\pi( \mathbf{b}^{(p)} | \mathbf{u}^{(p)} )
		\propto \exp\left( - \frac{1}{2} \norm{ H \mathbf{u}^{(p)} - \mathbf{b}^{(p)} }_{\Gamma}^2 \right).
\end{equation}
Observe that the data fidelity term in the traditional EnKF's analysis step~\cref{eq:EnKF_analysis2} is proportional to the negative log of~\cref{eq:likelihhod}.

\subsubsection{The conditional prior}
\label{subsub:cond_prior}

The joint prior $\pi( \mathbf{u}^{(p)}, \boldsymbol{\theta}^{(p)}|\what{\mathbf{u}}^{(p)} )$ in~\cref{eq:Bayes_rule} encodes our structural beliefs about $\mathbf{u}^{(p)}$.
Notably, we can decompose the joint prior as $\pi( \mathbf{u}^{(p)}, \boldsymbol{\theta}^{(p)} ) = \pi( \mathbf{u}^{(p)} | \boldsymbol{\theta}^{(p)},\what{\mathbf{u}}^{(p)}) \, \pi( \boldsymbol{\theta}^{(p)} )$ with conditional prior $\pi( \mathbf{u}^{(p)} | \boldsymbol{\theta}^{(p)},\what{\mathbf{u}}^{(p)})$ and hyper-prior $\pi( \boldsymbol{\theta}^{(p)} )$.%
\footnote{We assume $\boldsymbol{\theta}$ is independent of the forecast a priori, i.e., $\pi(\boldsymbol{\theta}^{(p)}|\what{\mathbf{u}}^{(p)})=\pi(\boldsymbol{\theta}^{(p)})$.}
We discuss the conditional prior below and the hyper-prior in \Cref{subsub:hyper_prior}.

The motivation for the {forecast} regularization term $\| \mathbf{u}^{(p)} - \what{\mathbf{u}}^{(p)} \|^2_{\what{C}}$ in the traditional EnKF's analysis step~\cref{eq:EnKF_analysis2} is that each particle approximately evolves according to the dynamics model~\cref{eq:dynamics_model}.
In a statistical setting, we model this as
\begin{equation}\label{eq:temporal_model}
	\mathbf{u}^{(p)} \sim \mathcal{N}( \what{\mathbf{u}}^{(p)}, \what{C} ).
\end{equation}
with forecast $\what{\mathbf{u}}^{(p)}$ as in~\cref{eq:EnKF_prediction} and sample prior covariance matrix $\what{C}$ as in~\cref{eq:EnKF_prior_covariance}, where we assume thus far that $\what{C}$ is invertible.
In practice, we replace $\what{C}$ by the localized and inflated sample covariance $\lambda\, (L \odot \what{C})$, where $L$ is a tapering matrix (see \Cref{app:details}) and $\lambda > 0$ is the forecast weight introduced in~\cref{eq:regEnKF_analysis}; a value $\lambda > 1$ corresponds to a multiplicative inflation of the forecast covariance within the analysis step.
To keep the notation light, we continue to write $\what{C}$ for this matrix in the remainder of this section.

The crux of our proposed GSBL-EnKF approach is to additionally incorporate the structural belief that the underlying system's states are piecewise smooth.
Within the GSBL framework~\cite{calvetti2020sparse,glaubitz2023generalized,glaubitz2026efficient}, we model $\mathbf{u}^{(p)}$ being piecewise smooth by $S \mathbf{u}^{(p)}$ being sparse for a sparsifying transform $S \in \R^{\nt \times \ns}$.
As mentioned before, $S$ can, for example, approximate a derivative operator.
To this end, we assume
\begin{equation}\label{eq:spatial_model}
	\pi(S \mathbf{u}^{(p)} | \boldsymbol{\theta}^{(p)}) = \mathcal{N}( \mathbf{0}, \Theta^{(p)} )
\end{equation}
with hyper-parameter vector $\boldsymbol{\theta}^{(p)} = [ \theta^{(p)}_1, \dots, \theta^{(p)}_\nt ]$ and diagonal matrix
$\Theta^{(p)} = \diag(\boldsymbol{\theta}^{(p)})$.

\begin{remark}\label{rem:motiv_cond_prior}
	The conditional Gaussian edge-preserving prior model~\eqref{eq:spatial_model} can be motivated by its asymptotic behavior~\cite{calvetti2007gaussian,glaubitz2023generalized}:
	Assume that $\theta^{(p)}_1=\dots=\theta^{(p)}_\nt$. Then, the model~\eqref{eq:spatial_model} favors values $\mathbf{u}^{(p)}$ for which the sparsified-state $S \mathbf{u}^{(p)}$ is close to zero, since such a state $\mathbf{u}^{(p)}$ has a higher probability.
	For instance, when $S \mathbf{u}^{(p)}$ corresponds to second differences of $\mathbf{u}^{(p)}$, i.e., $[ S \mathbf{u}^{(p)} ]_k = u^{(p)}_{k+1} - 2 u^{(p)}_{k} + u^{(p)}_{k-1}$ on an equidistant grid, then~\cref{eq:spatial_model} favors $\mathbf{u}^{(p)}$ to be close to linear.
	On the other hand, if a single hyperparameter $\theta^{(p)}_k$ is significantly larger than the others, a kink or jump of $\mathbf{u}^{(p)}$ near $x_k$ becomes more likely.
	In this way, the model~\cref{eq:spatial_model} promotes sparsity of $S \mathbf{u}^{(p)}$.
\end{remark}

We next formulate the conditional prior $\pi( \mathbf{u}^{(p)} | \boldsymbol{\theta}^{(p)},\what{\mathbf{u}}^{(p)})$ that mixes the above {forecast} and {edge-preserving} structural beliefs.
First, for $A\in\mathbb{R}^{m_A\times n_A}$ and $B\in\mathbb{R}^{m_B\times n_B}$, define the operator:
\begin{equation}
	\mathrm{diag}(A,B) \coloneqq
		\begin{bmatrix}
			A & 0 \\
			0 & B
		\end{bmatrix}
		\in\mathbb{R}^{(m_A + m_B) \times (n_A + n_B)}.
\end{equation}
Similarly, define $[A;B]\in\mathbb{R}^{(m_A + m_B) \times n}$ to be the vertical concatenation of the two matrices defined when $n_A = n_B = n$.
To get the appropriate conditional prior, we combine~\cref{eq:temporal_model,eq:spatial_model} to get the distribution:
\begin{equation}\label{eq:cond_prior_model}
	\pi(R \mathbf{u}^{(p)} \mid \boldsymbol{\theta}^{(p)},\what{\mathbf{u}}^{(p)}) = \mathcal{N}( \boldsymbol{\mu}^{(p)}, C_r^{(p)} ).
\end{equation}
Here, we denote transform $R = [I ; S] \in \R^{(\ns + \nt) \times \ns}$, mean $\boldsymbol{\mu}^{(p)} = [\what{\mathbf{u}}^{(p)}; \mathbf{0}] \in \R^{\ns + \nt}$, and covariance matrix $C_r^{(p)} = \diag( \what{C}, \Theta^{(p)} ) \in \R^{(\ns + \nt) \times (\ns + \nt)}$.
The density of the tilted Gaussian prior distribution is
\begin{equation}\label{eq:cond_prior_true_dens}
\begin{aligned}
\resizebox{.9\textwidth}{!}{$\displaystyle %
	\pi( \mathbf{u}^{(p)} | \boldsymbol{\theta}^{(p)},\what{\mathbf{u}}^{(p)})
		\propto
		\exp\left(
			- \frac{1}{2} \norm{ \mathbf{u}^{(p)} - \what{\mathbf{u}}^{(p)} }_{\what{C}}^2
			- \frac{1}{2} \norm{ S \mathbf{u}^{(p)} }_{\Theta^{(p)}}^2
			+ \frac{1}{2}\log\det\left(\what{C}^{-1} + S^\top \left(\Theta^{(p)}\right)^{-1} S\right)
		\right).
$}
\end{aligned}
\end{equation}
Observe that the negative log of density~\cref{eq:cond_prior_true_dens} can be interpreted as the sum of three terms. The forecast term, $\| \mathbf{u}^{(p)} - \what{\mathbf{u}}^{(p)} \|_{\what{C}}^2$, constrains the random variable $\mathbf{u}^{(p)}$ to lie near the simulations.
The edge-preserving regularization term, $\| S \mathbf{u}^{(p)} \|_{\Theta^{(p)}}^2$, includes a spatially-varying parameter $\boldsymbol{\theta}^{(p)}$, which aligns with the discussion of regularization presented in \Cref{sec:basic_idea}.
Finally, we note that~\cref{eq:cond_prior_true_dens} defines the conditional prior as the product of the forecast factor~\cref{eq:temporal_model} and the regularity factor~\cref{eq:spatial_model}: Since $R = [I; S]$ is not square, \cref{eq:cond_prior_model} defines a degenerate conditional Gaussian on $\R^{\ns+\nt}$.
Unfortunately, the exact log-normalization of the product in $\mathbf{u}^{(p)}$, $\log\det(\what{C}^{-1} + S^\top (\Theta^{(p)})^{-1} S)$, couples the hyperparameters in a computationally non-trivial way.
Following the GSBL literature~\cite{calvetti2020sparse,glaubitz2023generalized,xiao2023sequential,lindbloom2024generalized}, however, we drop this coupling term in favor of a more tractable choice.
Specifically, we replace $\log\det\left(\what{C}^{-1} + S^\top \left(\Theta^{(p)}\right)^{-1} S\right)$ by $\log\det({\Theta^{(p)}}^{-1}) = -\sum_{k=1}^{n_s}\log\theta^{(p)}_k$, which is what renders the hyperparameter updates in \Cref{subsub:theta_update} explicit and decoupled.
The density we opt to sample from is thus given as:

\begin{equation}\label{eq:cond_prior}
\begin{aligned}
	\tilde{\pi}( \mathbf{u}^{(p)} | \boldsymbol{\theta}^{(p)},\what{\mathbf{u}}^{(p)})
		\propto \exp\left(
			- \frac{1}{2} \norm{ \mathbf{u}^{(p)} - \what{\mathbf{u}}^{(p)} }_{\what{C}}^2
			- \frac{1}{2} \norm{ S \mathbf{u}^{(p)} }_{\Theta^{(p)}}^2 - \frac{1}{2}\sum_{k=1}^{n_s}\log\theta^{(p)}
		\right).
\end{aligned}
\end{equation}

\subsubsection{The hyper-prior}
\label{subsub:hyper_prior}

From the discussion in \Cref{rem:motiv_cond_prior}, it is evident that the variance hyper-parameters $\theta^{(p)}_1, \dots, \theta^{(p)}_\nt$ must be allowed to have values that vary across the entries of the conditional Gaussian prior model~\cref{eq:spatial_model} to promote sparsity for $S \mathbf{u}^{(p)}$. Specifically, we promote the values $\theta^{(p)}_k$ for $k=1,\ldots,\nt$ being nearly zero for smooth segments, while allowing for a few occasional outliers corresponding to jump discontinuities.
This can be achieved by treating $\theta^{(p)}_1, \dots, \theta^{(p)}_\nt$ as independent random variables under the prior, where the prior density:
\begin{itemize}
	\item[(i)]
	is supported on the positive real numbers;

	\item[(ii)] concentrates its mass at small values of $\theta$; and

	\item[(iii)]
	has a heavy-tailed decay.

\end{itemize}
Condition (i) ensures that the covariance matrix $\Theta^{(p)} = \diag(\boldsymbol{\theta}^{(p)})$ in~\cref{eq:spatial_model} is positive definite. Condition (ii) promotes near-zero values of $[S\mathbf{u}^{(p)}]_k$ due to rare deviations from zero.
Condition (iii) allows for occasional outliers, i.e., there can be a few values of $k$ such that entry $\theta^{(p)}_k$ has a significantly larger value than most other entries of the ensemble member $\boldsymbol{\theta}^{(p)}$. The corresponding entries of $S \mathbf{u}^{(p)}$ are then more likely to differ from zero.
A common choice satisfying the above desiderata is the generalized gamma distribution~\cite{calvetti2020sparse}, which has the density:
\begin{equation}
	\mathcal{GG}( \theta | r, \beta, \vartheta )
		= \frac{|r|}{\Gamma(\beta) \vartheta} \left[ \frac{\theta}{\vartheta} \right]^{r\beta-1} \exp\left( - \left[ \frac{\theta}{\vartheta} \right]^r \right),
\end{equation}
where the parameters are the rate $r \in \R \setminus \{0\}$, the shape $\beta > 0$, and the scale $\vartheta > 0$.
We assume that the hyper-parameters independently and identically follow the generalized gamma. That is, the tensorization of hyper-prior $\pi(\theta^{(p)}_k) = \mathcal{GG}(r, \beta, \vartheta)$ across $k=1,\ldots,\nt$ yields the joint density:
\begin{equation}\label{eq:hyper_prior}
	\pi( \boldsymbol{\theta}^{(p)} )
		\propto \exp\left( - \sum_{k=1}^\nt \left[ \frac{\theta^{(p)}_k}{\vartheta} \right]^r  + (r\beta - 1)\sum_{k=1}^{\nt} \log\theta^{(p)}_k \right).
\end{equation}
Using a generalized gamma hyper-prior allows for particularly efficient inference~\cite{calvetti2020sparsity,lindbloom2024generalized}, which we will leverage in \Cref{sub:analysis_step}.

\subsubsection{The posterior}
\label{subsub:posterior}

We are now positioned to formulate the posterior~\cref{eq:Bayes_rule} underpinning the proposed GSBL-EnKF approach.
To this end, substituting~\cref{eq:likelihhod,eq:cond_prior,eq:hyper_prior} into~\cref{eq:Bayes_rule} yields
\begin{equation}\label{eq:posterior}
\begin{aligned}
	\tilde{\pi}( \mathbf{u}^{(p)}, \boldsymbol{\theta}^{(p)} | \mathbf{b}^{(p)},\what{\mathbf{u}}^{(p)})
		\propto \exp\Bigg(
			& - \frac{1}{2} \norm{ H \mathbf{u}^{(p)} - \mathbf{b}^{(p)} }_{\Gamma}^2
			- \frac{1}{2} \norm{ \mathbf{u}^{(p)} - \what{\mathbf{u}}^{(p)} }_{\what{C}}^2
			- \frac{1}{2} \norm{ S \mathbf{u}^{(p)} }_{\Theta^{(p)}}^2 \\
			& - \sum_{k=1}^\nt \left[ \theta^{(p)}_k / \vartheta \right]^r
			+ \tau \sum_{k=1}^\nt \log \theta^{(p)}_k
		\Bigg)
\end{aligned}
\end{equation}
for the posterior density of the $p$-th particle, where we define the constant $\tau\coloneqq r \beta - 3/2 $ and denote it $\tilde{\pi}$ due to the approximation to the normalization constant made in~\cref{eq:cond_prior}.

\subsection{The GSBL-EnKF analysis step}
\label{sub:analysis_step}

We now formulate the corresponding GSBL-EnKF analysis step as the MAP estimate of the joint posterior~\cref{eq:posterior}.
This allows us to interpret the GSBL-EnKF's analysis step naturally as a regularized version of the traditional EnKF's analysis step~\cref{eq:EnKF_analysis2}.

Recall that $( \mathbf{u}^{(1:P)}_{\text{MAP}}, \boldsymbol{\theta}^{(1:P)}_{\text{MAP}} )$ is the MAP estimate of the joint distribution $\tilde\pi( \mathbf{u}^{(1:P)}, \boldsymbol{\theta}^{(1:P)} | \mathbf{b}^{(1:P)}, \what{\mathbf{u}}^{(1:P)} )$ if it maximizes the density of this posterior.
Equivalently, yet often computationally more robust, we find $( \mathbf{u}^{(1:P)}_{\text{MAP}}, \boldsymbol{\theta}^{(1:P)}_{\text{MAP}} )$ by minimizing the negative log-posterior:
\begin{equation}\label{eq:MAP_estimate}
	( \mathbf{u}^{(1:P)}_{\text{MAP}}, \boldsymbol{\theta}^{(1:P)}_{\text{MAP}} )
		= \argmin_{ \mathbf{u}^{(1:P)}, \boldsymbol{\theta}^{(1:P)} } \left\{ \mathcal{J}_P( \mathbf{u}^{(1:P)}, \boldsymbol{\theta}^{(1:P)} ) \right\},
\end{equation}
with the objective function $\mathcal{J}_P( \mathbf{u}^{(1:P)}, \boldsymbol{\theta}^{(1:P)} ) = - \log\tilde{\pi}( \mathbf{u}^{(1:P)}, \boldsymbol{\theta}^{(1:P)} | \mathbf{b}^{(1:P)} ) + C$, given by
\begin{gather}
    \mathcal{J}_{P}(\mathbf{u}^{(1:P)}, \boldsymbol{\theta}^{(1:P)}) := \sum_{p=1}^P \mathcal{J}(\mathbf{u}^{(p)},\boldsymbol{\theta}^{(p)};\mathbf{b}^{(p)},\what{\mathbf{u}}^{(p)}),\\[-0.5em]
    \mathcal{J}(\mathbf{u},\boldsymbol{\theta};\mathbf{b},\what{\mathbf{u}}) := \frac{1}{2} \norm{H\mathbf{u} - \mathbf{b}}^2_{\Gamma} + \frac{1}{2}\norm{\mathbf{u} - \what{\mathbf{u}}}^2_{\what{C}} + \frac{1}{2}\norm{S\mathbf{u}}^2_{\Theta} + \sum_{k=1}^{\nt}\left( \left[\theta_k / \vartheta\right]^r - [r\beta - 3/2]\log\theta_k \right), \nonumber
\end{gather}
where the constant $C$ is independent of $\mathbf{u}^{(1:P)}$ and $\boldsymbol{\theta}^{(1:P)}$.

We use a block-coordinate descent approach~\cite{wright2015coordinate,beck2017first} to approximate the MAP estimate~\cref{eq:MAP_estimate} efficiently.
A prevalent block-coordinate descent method for the joint posterior in the GSBL context is the so-called iterative alternating sequential (IAS) algorithm~\cite{calvetti2020sparse,calvetti2020sparsity,glaubitz2023leveraging,lindbloom2024generalized}.
The IAS algorithm computes the minimizer of the objective function $\mathcal{J}$ by alternating between (i) minimizing $\mathcal{J}$ w.r.t.\ $\mathbf{u}^{(1:P)}$ for fixed $\boldsymbol{\theta}^{(1:P)}$ and (ii) minimizing $\mathcal{J}$ w.r.t.\ $\boldsymbol{\theta}^{(1:P)}$ for fixed $\mathbf{u}^{(1:P)}$.
To ensure clarity, we will use an underline to designate a variable as fixed.
Then, given an initial guess $\underline{\boldsymbol{\theta}}^{(1:P)}$ for the parameters, the IAS algorithm proceeds through a sequence of updates of the form
\begin{equation}\label{eq:IAS}
	\underline{\mathbf{u}}^{(1:P)}
		= \argmin_{\mathbf{u}^{(1:P)}} \left\{ \mathcal{J}_P(\mathbf{u}^{(1:P)}, \underline{\boldsymbol{\theta}}^{(1:P)}) \right\}, \quad
	\underline{\boldsymbol{\theta}}^{(1:P)}
		= \argmin_{\boldsymbol{\theta}^{(1:P)}} \left\{ \mathcal{J}_P(\underline{\mathbf{u}}^{(1:P)}, \boldsymbol{\theta}^{(1:P)}) \right\},
\end{equation}
until a convergence criterion is met.
The IAS algorithm is an attractive choice because the two subproblems~\cref{eq:IAS} are significantly easier to solve than the joint problem~\cref{eq:MAP_estimate}. We note that the decomposition of $\mathcal{J}_P$ as the sum of terms $\mathcal{J}$ allows us to update each state realization $\mathbf{u}^{(p)}$ given fixed parameters $\underline{\boldsymbol{\theta}}^{(p)}$ in parallel across the ensemble members, $p=1,\ldots,P$.
The same holds for updating parameters $\boldsymbol{\theta}^{(p)}$ given fixed state $\underline{\mathbf{u}}^{(p)}$ in parallel.

\subsection{Updating the states \texorpdfstring{$\mathbf{u}^{(1:P)}$}{}}
\label{subsub:u_update}

Updating $\mathbf{u}^{(1:P)}$ given $\underline{\boldsymbol{\theta}}^{(1:P)}$ reduces to solving the decoupled quadratic minimization problems:
\begin{equation}\label{eq:u_update1}
	\mathbf{u}^{(p)}
		= \argmin_{ \mathbf{u} \in \R^\ns } \left\{
			\norm{ H \mathbf{u} - \mathbf{b}^{(p)} }_{\Gamma}^2
			+ \norm{ \mathbf{u} - \what{\mathbf{u}}^{(p)} }_{\what{C}}^2
			+ \norm{ S \mathbf{u} }_{\underline{\Theta}^{(p)}}^2
		\right\},\quad \underline{\Theta}^{(p)} \coloneqq \mathrm{diag}(\underline{\boldsymbol{\theta}}^{(p)}),
\end{equation}
for $p=1,\dots,P$, where we have ignored terms independent of $\mathbf{u}^{(1:P)}$.
Solving~\cref{eq:u_update1} is equivalent to solving for the least squares solutions of the overdetermined linear systems
\begin{equation}\label{eq:u_update2}
    \mathbf{u}^{(p)} = \argmin_{\mathbf{u}\in\mathbb{R}^{\ns}}\left\|\begin{bmatrix}
    	\Gamma^{-1/2} H \\
		\what{C}^{-1/2} I \\
		\left(\underline{\Theta}^{(p)}\right)^{-1/2} S
	\end{bmatrix}
	\mathbf{u}
     -
     \begin{bmatrix}
     	\Gamma^{-1/2} \mathbf{b}^{(p)} \\
		\what{C}^{-1/2} \what{\mathbf{u}}^{(p)} \\
		\mathbf{0}
	\end{bmatrix}\right\|_2^2 ,\quad p=1,\dots,P,
\end{equation}%
which, in turn, is equivalent to solving the linear system:
\begin{equation}\label{eq:u_update3}
    \underbrace{\left( H^\top \Gamma^{-1} H
    		+ \what{C}^{-1}
		+ S^\top \left(\underline{\Theta}^{(p)}\right)^{-1} S \right)
	}_{:= \Sigma^{(p)}} \mathbf{u}^{(p)}
    	= H^\top \Gamma^{-1} \mathbf{b}^{(p)} + \what{C}^{-1} \what{\mathbf{u}}^{(p)}, \quad p=1,\dots,P.
\end{equation}
Notably, if the coefficient matrix $\Sigma^{(p)}$ on the left-hand side of~\cref{eq:u_update3} is symmetric positive-definite (SPD), each of~\cref{eq:u_update1,eq:u_update2,eq:u_update3} shares the same unique solution.

\subsection{Connection to the traditional Kalman update}
\label{sub:Kalman}
Unfortunately, the given updates~\cref{eq:u_update1,eq:u_update2,eq:u_update3} require the (Cholesky decomposition of the) inverse covariance matrices $\Gamma^{-1}$ and $\what{C}^{-1}$.
Therefore, we next formulate the $\mathbf{u}^{(1:P)}$-update~\cref{eq:u_update1,eq:u_update2,eq:u_update3} as a regularized version of the EnKF Kalman update~\cref{eq:EnKF_analysis1}.

Computing matrix inverses, however, is often computationally impractical or even impossible since $\what{C}$ can be low-rank (e.g., when $P < \ns$) or ill-conditioned.
The Kalman-like update we provide below sidesteps this issue by only requiring the matrices $\Gamma$ and $\what{C}$, rather than their inverses.
In our experiments, $\what{C}$ is a tapered sample covariance $L \odot \what{C}$, which is SPD whenever the tapering matrix $L$ is SPD and all ensemble variances are positive (Schur product theorem).
The following lemma is thus stated for SPD $\what{C}$.

\begin{lemma}[Kalman formulation of the $\mathbf{u}^{(1:P)}$-update]\label{lem:Kalman}
    Assume that $\what{C}$ is SPD, and let $K_r^{(p)} = \what{C} H_r^\top (H_r \what{C} H_r^\top + \Gamma_r^{(p)} )^{-1}$ be the \emph{regularized Kalman gain matrices} with $H_r = [H; S] \in \R^{(\no+\nt) \times \ns}$ and $\Gamma_r^{(p)} = \diag( \Gamma, \underline{\Theta}^{(p)} ) \in \R^{(\no+\nt) \times (\no+\nt)}$.
	The $\mathbf{u}^{(1:P)}$-updates~\cref{eq:u_update1,eq:u_update2,eq:u_update3} can be written:
    \begin{equation}\label{eq:x_update_Kalman}
        \mathbf{u}^{(p)}
            = \what{\mathbf{u}}^{(p)} + K_r^{(p)} \left(
            \begin{bmatrix}
            		\mathbf{b}^{(p)} \\ \mathbf{0}
            \end{bmatrix}
            - H_r \what{\mathbf{u}}^{(p)}
            \right),
            \quad p=1,\dots,P.
    \end{equation}
\end{lemma}

The proof of \Cref{lem:Kalman} follows from arguments used in the traditional Kalman update~\cref{eq:EnKF_analysis1}; we provide a proof in \Cref{app:proofs}.

\begin{remark}
	Compared to the EnKF, the regularized update~\cref{eq:x_update_Kalman} is more expensive in two respects.
	First, the dimension of the linear system grows from $\no$ to $\no + \nt$.
	Second, since $\Gamma_r^{(p)}$ depends on $p$, a separate system has to be solved for every ensemble member and every IAS iteration, i.e., $(\nIAS + 1) P$ solves per analysis step instead of one factorization shared by all members.
	In large-scale simulations, the Kalman gain $K_r^{(p)}$ would therefore not be formed explicitly but applied using Krylov methods that only require matrix-vector products with $H_r \what{C} H_r^\top + \Gamma_r^{(p)}$.
	Alternatively, when $\what{C}$ is SPD, the information form~\cref{eq:u_update3} can be used, which is of dimension $\ns$ and shares a single factorization of $\what{C}$ across all members and iterations.
	This added complexity of the algorithm can be reduced, however, by following~\cite{glaubitz2023leveraging} and replacing the ensemble of hyperparameters $\boldsymbol{\theta}^{(1:P)}$ with a single set of parameters $\boldsymbol{\theta}$ shared across all ensemble members $\mathbf{u}^{(1:P)}$, thus keeping the Kalman gain independent of ensemble index $p$.
 	While this alternative choice for sharing parameters is more computationally convenient, it eliminates a key component of modeling epistemic uncertainty: the assumption that one can choose regularization hyperparameters to be constant across the entire ensemble is tantamount to suggesting that there is little-to-no uncertainty about the location of a state's discontinuity.
 	We therefore choose to model each ensemble member with its own hyperparameter to represent uncertainty about the discontinuity location.
\end{remark}

\subsection{Updating the hyper-parameters\texorpdfstring{ $\boldsymbol{\theta}$}{}}
\label{subsub:theta_update}
Updating a hyper-parameter vector $\boldsymbol{\theta}^{(p)} = [ \theta_1, \dots, \theta_\nt ]$ given fixed $\underline{\mathbf{u}}^{(p)}$ requires solving the decoupled minimization problems
\begin{equation}\label{eq:theta_update1}
	\theta_{k}^{(p)}
		= \argmin_{ \theta_k > 0 } \left\{
			\frac{[ S \underline{\mathbf{u}}^{(p)} ]_k^2}{2\theta_k}  + \left( \frac{\theta_k}{\vartheta} \right)^{r} - \tau \log \theta_k
		\right\}, \quad k=1,\dots,\nt,
\end{equation}
where we have ignored all terms that do not depend on $\boldsymbol{\theta}$.
Furthermore, recall that $[S\mathbf{u}]_k$ denotes the $k$-th entry of the vector $S\mathbf{u} \in \R^\nt$ and recall parameter $\tau = r \beta - 3/2$.
Differentiating the objective function in~\cref{eq:theta_update1} w.r.t.\ $\theta_k$ and setting this derivative to zero yields
\begin{equation}\label{eq:theta_update2}
	0 = - \frac{[ S \underline{\mathbf{u}}^{(p)} ]_k^2}{2\theta_k^2}
		+ \theta_k^{r-1} \left( \frac{r}{\vartheta^{r}} \right)
		- \frac{\tau}{\theta_k}, \quad k=1,\dots,\nt.
\end{equation}
Following~\cite{calvetti2020sparse,calvetti2020sparsity}, one can show that~\cref{eq:theta_update2} produces the update rule:
\begin{equation}\label{eq:theta_update3}
	\theta_k^{(p)} =
		\vartheta \cdot \varphi\left( \frac{ \big|\, [ S \underline{\mathbf{u}}^{(p)} ]_k \big| }{ \sqrt{\vartheta} } \right), \quad k=1,\dots,\nt,
\end{equation}
where $\varphi$ is the solution to the following scalar initial value problem:
\begin{equation}\label{eq:eq:theta_update_IVP}
	\varphi'(t) = \cfrac{2 t \varphi(t)}{ 2 r^2 \varphi(t)^{r + 1} + t^2 }, \quad
	\varphi(0) = \left( \frac{\tau}{r} \right)^{1/r},
\end{equation}
assuming that either (i) $r < 0$ or (ii) $r > 3/(2\beta)$.
The updates in~\cref{eq:theta_update3} can be efficiently calculated by numerically solving~\cref{eq:eq:theta_update_IVP} only once a priori, i.e., without any evaluation of $\mathbf{u}$ or $\boldsymbol{\theta}$, over a sufficiently long time horizon, then evaluating it at the appropriate points when performing the GSBL-EnKF.
Moreover, for $r = \pm1$, the condition~\cref{eq:theta_update2} admits a simple explicit solution formula~\cite{calvetti2020sparse,calvetti2020sparsity,glaubitz2023leveraging,lindbloom2024generalized}.

\subsection{Algorithm summary and discussion}
\label{sub:discussion}

We summarize the implementation of the proposed GSBL-EnKF approach in \Cref{alg:GSBL_EnKF}.
\begin{algorithm}[ht!]
\caption{GSBL-EnKF method}\label{alg:GSBL_EnKF}
\begin{algorithmic}
	\State{\textbf{Input:}} Initial ensemble $\{ \mathbf{u}_{0}^{(p)} \}_{p=1}^P$
	\For{$j=1,\dots,J$}
		\State{\textbf{Prediction step.}}
			\State Compute $\what{\mathbf{u}}_{j}^{(p)} = \Psi( \mathbf{u}_{j-1}^{(p)} )$ for $p=1,\dots,P$
			\State Compute the prior mean $\what{\mathbf{m}}_{j} = \frac{1}{P} \sum_{p=1}^P \what{\mathbf{u}}_{j}^{(p)}$
			\State Compute the prior covariance $\what{C}_{j} = \frac{1}{P-1} \sum_{p=1}^P ( \what{\mathbf{u}}_{j}^{(p)} - \what{\mathbf{m}}_{j} ) ( \what{\mathbf{u}}_{j}^{(p)} - \what{\mathbf{m}}_{j} )^\top$
		\State{\textbf{Analysis step.}}
			\State Generate perturbed observations $\mathbf{b}_{j}^{(p)} = \mathbf{b}_{j} + \boldsymbol{\eta}_{j}^{(p)}$, $\boldsymbol{\eta}_{j}^{(p)} \sim \mathcal{N}(\mathbf{0}, \Gamma)$ for $p=1,\dots,P$
			\State Initialize hyper-parameter vectors $\boldsymbol{\theta}_{j}^{(p)}= \boldsymbol{\theta}_0$ (we use $\boldsymbol{\theta}_0 = \mathbf{1}$; see \Cref{app:details})
			\State Replace $\what{C}_j$ by the localized and scaled covariance $\lambda\,(L \odot \what{C}_j)$ (\Cref{app:details})
			\For{$i = 1,\dots,\nIAS$}
				\State Update particles $\mathbf{u}_{j}^{(1:P)}$ according to~\cref{eq:x_update_Kalman}
				\State Update the hyper-parameter vectors $\boldsymbol{\theta}_{j}^{(p)}$ according to~\cref{eq:theta_update3}
			\EndFor
            \State Finalize particle positions $\mathbf{u}_{j}^{(1:P)}$ according to~\cref{eq:x_update_Kalman} and perform any inflation
	\EndFor
	\State{\textbf{Output:}} Ensembles $\{ \mathbf{u}_j^{(p)} \}_{p=1}^P$ and hyper-parameter vectors $\{ \boldsymbol{\theta}_{j}^{(p)} \}_{p=1}^P$ for $j=0,\dots,J$
\end{algorithmic}
\end{algorithm}

Contrasting the proposed GSBL-EnKF framework in \Cref{alg:GSBL_EnKF} with the traditional EnKF, we note that the standard, unregularized Kalman update $\mathbf{u}_{j}^{(p)} = ( I - K_{j} H ) \what{\mathbf{u}}_{j}^{(k)} + K_{j} \mathbf{b}_{j}^{(p)}$ has been replaced by an iterative scheme. %
This IAS scheme alternates between the GSBL-regularized Kalman update~\cref{eq:x_update_Kalman} and updating the hyper-parameter vectors $\boldsymbol{\theta}_{j}^{(p)}$ for a fixed number of iterations $\nIAS \in \N$.
We observe that this iterative procedure converges rapidly; therefore, we keep $\nIAS$ small to avoid significantly increasing computational cost.
In our implementation, we typically choose $\nIAS = 2$.
Alternatively, a different convergence criterion could be used instead of a fixed iteration count.

\section{Numerical tests}
\label{sec:tests}

We use nodal discontinuous Galerkin solvers~\cite{hesthaven2007nodal,gassner2016split} for semi-discretizing each hyperbolic conservation law and the subcell shock-capturing strategy described in~\cite{hennemann2021provably}.
For time integration, we use a strong-stability preserving (SSP) Runge--Kutta (RK) method with adaptive time stepping~\cite{kraaijevanger1991contractivity,conde2018embedded,ranocha2022optimized};
the implementations are as given in~\cite{ranocha2022adaptive,rackauckas2017differentialequations}.
All assimilation is performed at the discretization's quadrature points, which we also use to estimate error via the time-averaged root-mean-squared error (RMSE) and the continuous ranked probability score (CRPS).
See \Cref{app:details} for discretization details.

\subsection{The linear advection equation}
\label{sub:tests_linear}

\begin{figure}[tb]
    \centering
    \includegraphics[width=0.5\linewidth,clip,trim={0.1cm 0.1cm 0.1cm 0.1cm}]{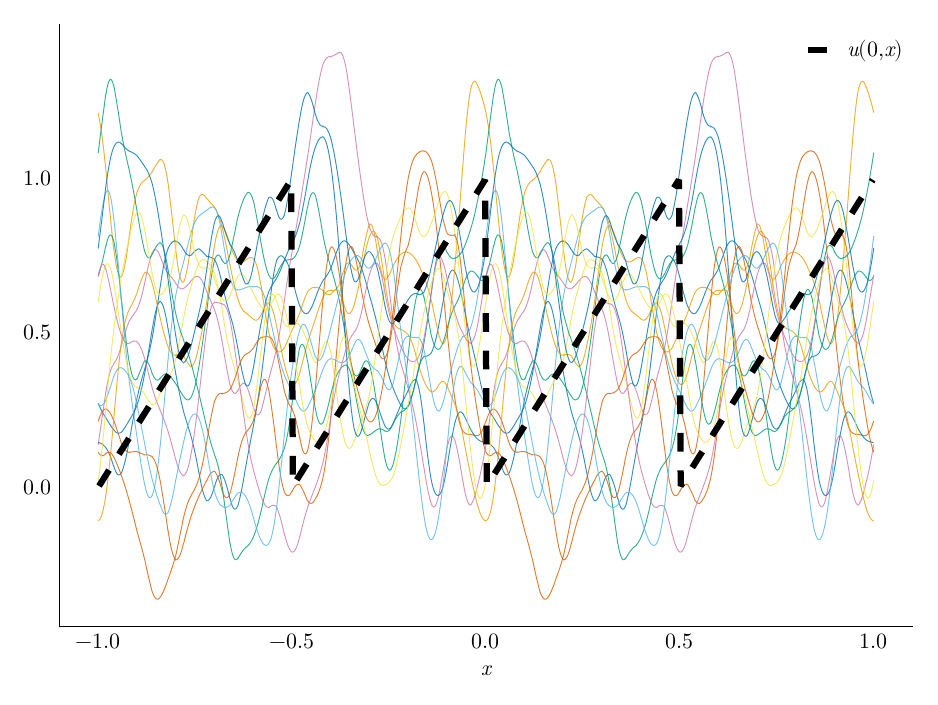}
    \caption{Initialization of the ensemble for the linear advection equation of \Cref{sub:tests_linear}.
    The first 15 ensemble members are shown alongside the data-generating initial condition.}
    \label{fig:advection_initial}
\end{figure}

Consider the linear advection equation
\begin{equation}\label{eq:linear}
    \partial_t u(t,x) + \alpha \partial_x u(t,x) = 0
\end{equation}
on $\Omega = [-1,1]$ with constant velocity $\alpha > 0$, periodic boundary conditions $u(t,0) = u(t,1)$, and sawtooth initial data $u(0,x) = \frac{x+1}{2}\;\mathrm{mod}\;\frac{1}{\tau},$ where ``$\mathrm{mod}$'' is the usual modulo operation and $\tau$ is an integer number of teeth.
We set $\alpha = 0.1$ and $\tau = 4$, recovering the exact solution at time $t$ as $u(t,x) = u(0, x-\alpha t)$.
That is, the solution remains piecewise linear, exposing the behavior of the proposed GSBL-EnKF approach over time.
We initialize the ensemble with smooth, oscillatory ensemble members, shown in \Cref{fig:advection_initial}. By initializing without knowledge of the number, locations, or amplitudes of the state discontinuities, we aim to evaluate the robustness of the data assimilation methods.
More information can be found in \Cref{app:details}.

\begin{figure}[tb]
    \centering
    \includegraphics[width=0.8\linewidth,clip,trim={0.3cm 0.9cm 0.5cm 0.5cm}]{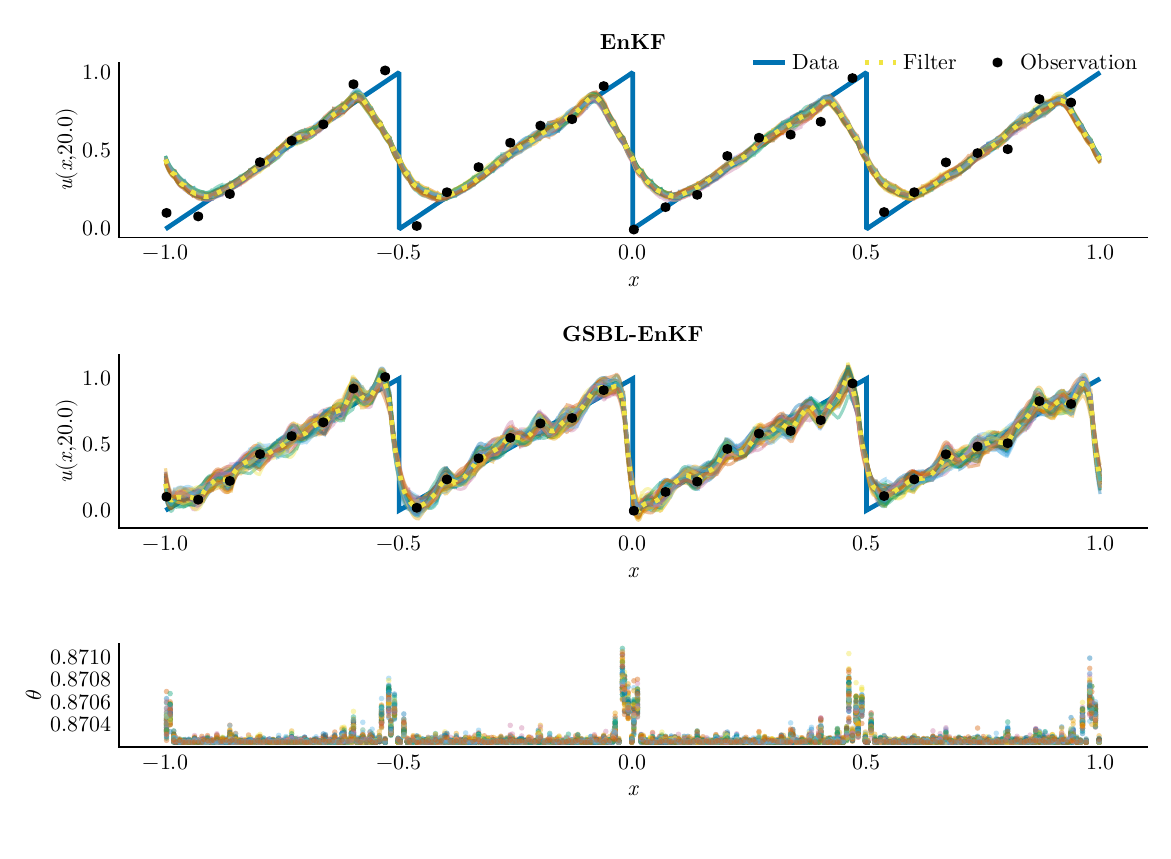}
    \caption{
        Time profiles of the EnKF vs.\ our proposed GSBL-EnKF for {the linear advection equation in} \Cref{sub:tests_linear} at time $t=10$. We use $P=40$ ensemble members and observation spacing $\Delta y\approx 0.07$.
        Our proposed GSBL-EnKF better preserves sharp jump discontinuities that the EnKF smears out.
    }
    \label{fig:advection_profiles}
\end{figure}

\begin{figure}[tb]
    \centering
    \includegraphics[width=0.99\linewidth,clip,trim={0.1cm 0.45cm 0.1cm 0.1cm}]{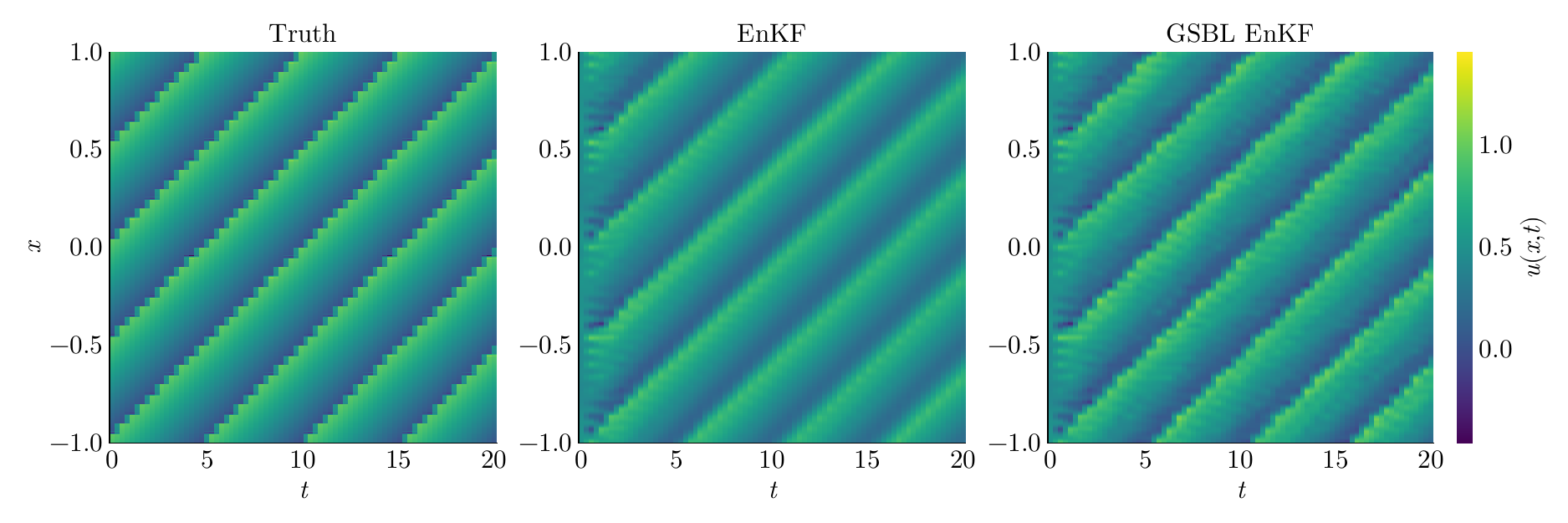}
    \caption{
        Comparing the data-generating flow to the ensemble means of {the EnKF and our proposed GSBL-EnKF for the linear advection simulation} in \Cref{fig:advection_profiles}.
        Our proposed GSBL-EnKF better preserves sharp jump discontinuities that the EnKF smears out.
        The staircase-like artifacts come from the discrete-time nature of our observations and do not reflect any quality of the methods.
        }
    \label{fig:advection_heatmap}
\end{figure}

\Cref{fig:advection_profiles} illustrates a single time profile of the ensemble members and the GSBL ensemble's found hyperparameters.
We observe that the hyperparameters $\theta^{(p)}_k$ peak at the jump locations of their respective ensemble members, with the spread of these peaks across members reflecting uncertainty about the jump locations.
The variation of $\theta^{(p)}_k$ is small in absolute terms because the transform $S$ is element-local (\Cref{rem:sparsifying_transform}) and the sawtooth is linear within elements, except in elements containing a jump.
\Cref{fig:advection_heatmap} shows the space-time plot of the average ensemble member.
We see that our GSBL-EnKF can achieve a more accurate solution prediction faster while ensuring that the shock is not smeared.

\begin{figure}[tb]
    \centering
    \includegraphics[width=0.8\linewidth,clip,trim={0.1cm 0.15cm 0.1cm 0.1cm}]{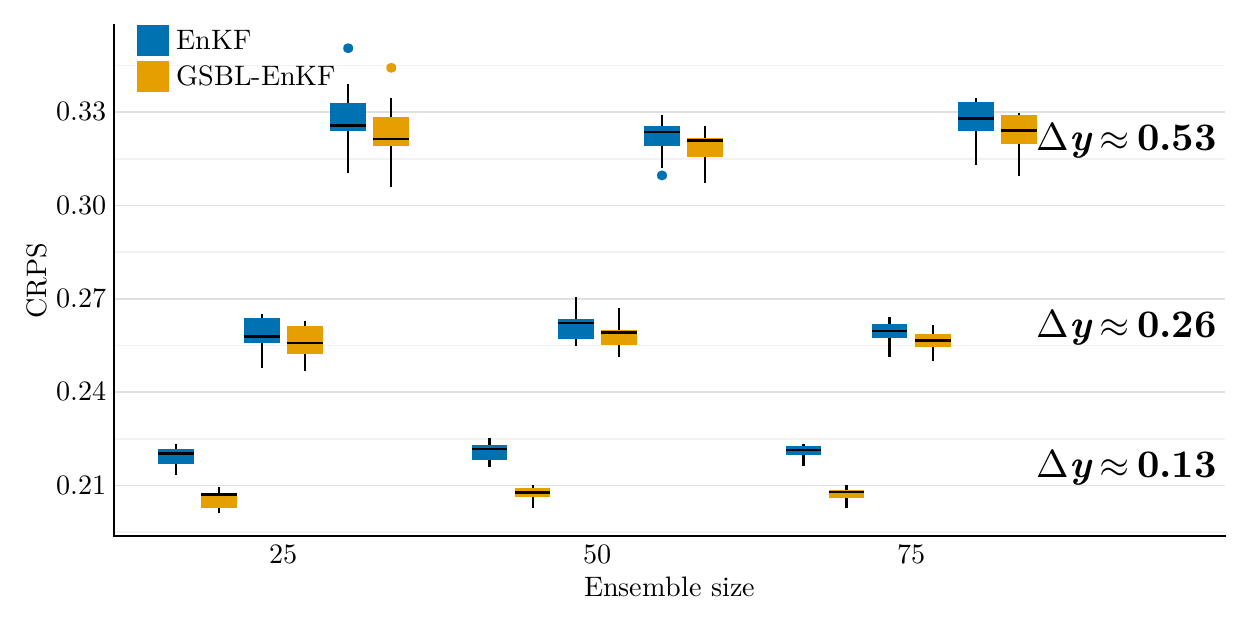}
    \caption{{The CRPS of the EnKF and our proposed GSBL-EnKF for different grid spacings and ensemble sizes across independent trials of the linear advection equation in \Cref{sub:tests_linear}.}}
    \label{fig:advection_convergence}
\end{figure}

We next quantitatively investigate the effect of changing several hyperparameters in the filtering problem for the linear advection equation in \Cref{fig:advection_convergence}.
Specifically, we report how the CRPS is affected by the ensemble size and the observation scaling. Notably, while our proposed GSBL-EnKF yields smaller CRPS values in all cases, the improvement is most notable for dense observation grids.
Across all ninety initializations---ten independent trials for each pair of observation grid and ensemble size---we observe that, for each initial ensemble, the GSBL-EnKF improves over the standard EnKF in both RMSE and CRPS.

\subsection{The inviscid Burgers equation}
\label{sub:tests_Burgers}

\begin{figure}[tb]
    \centering
    \includegraphics[width=0.5\linewidth,clip,trim={0.1cm 0.4cm 0.1cm 0.1cm}]{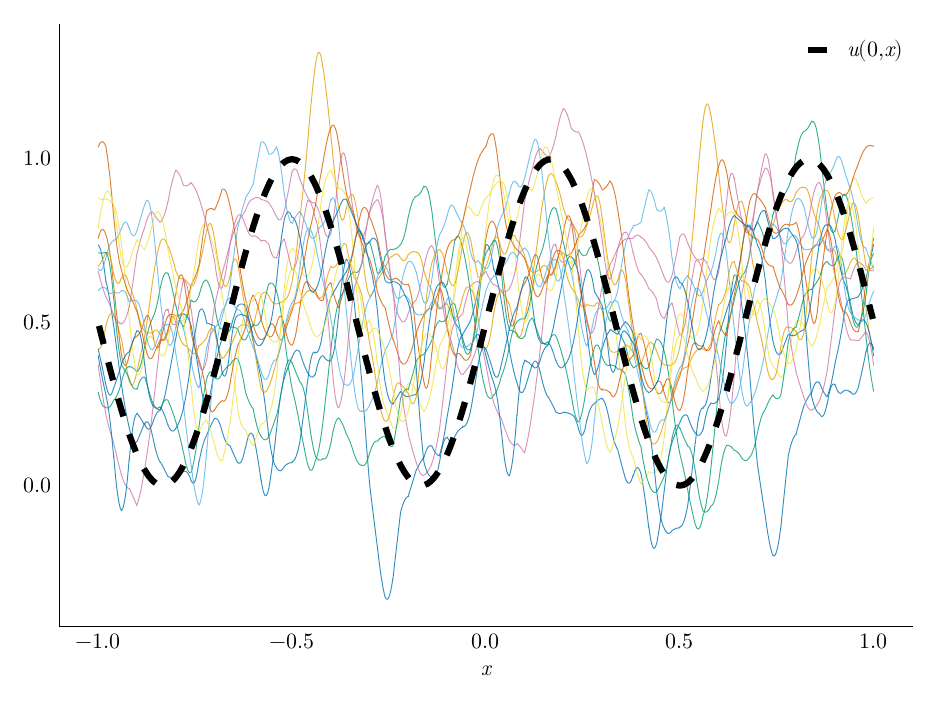}
    \caption{Initial ensemble for the inviscid Burgers equation in \Cref{sub:tests_Burgers}.}
    \label{fig:init_Burgers}
\end{figure}

Consider the inviscid Burgers equation
\begin{equation}\label{eq:Burgers}
    \partial_t u(t,x) + \partial_x \left( \frac{u(t,x)^2}{2} \right) = 0
\end{equation}
on $\Omega = [-1,1]$ with periodic boundary conditions $u(t,-1) = u(t,1)$, and smooth initial data $u(0,x) = 0.5 + 0.5\sin( 3 \pi x )$.
The Burgers equation is well known for developing piecewise-smooth solutions with jump discontinuities, even for smooth initial states and in finite time.
We thus use this more challenging test case to gauge the performance of our methods in systems that start smooth but develop shocks midway through the simulation.
We again initialize using smooth, oscillatory ensemble members, shown in \Cref{fig:init_Burgers} and described in \Cref{app:details}.
Unlike the advection example in \Cref{sub:tests_linear}, the Burgers equation becomes dissipative once a shock has formed; any information that travels into the shock is destroyed irreversibly.
In addition, we make the dynamics stochastic by adding white noise with standard deviation $\sigma_x$ to the state to study the effect of stochastic dynamics on the filtering methods.
In particular, the resulting inference problem can then no longer be solved by inferring the initial condition alone (cf. \Cref{sub:problems}).

\begin{figure}[ht!]
    \centering
    \includegraphics[width=0.8\linewidth,clip,trim={0.1cm 1cm 0.1cm 0.1cm}]{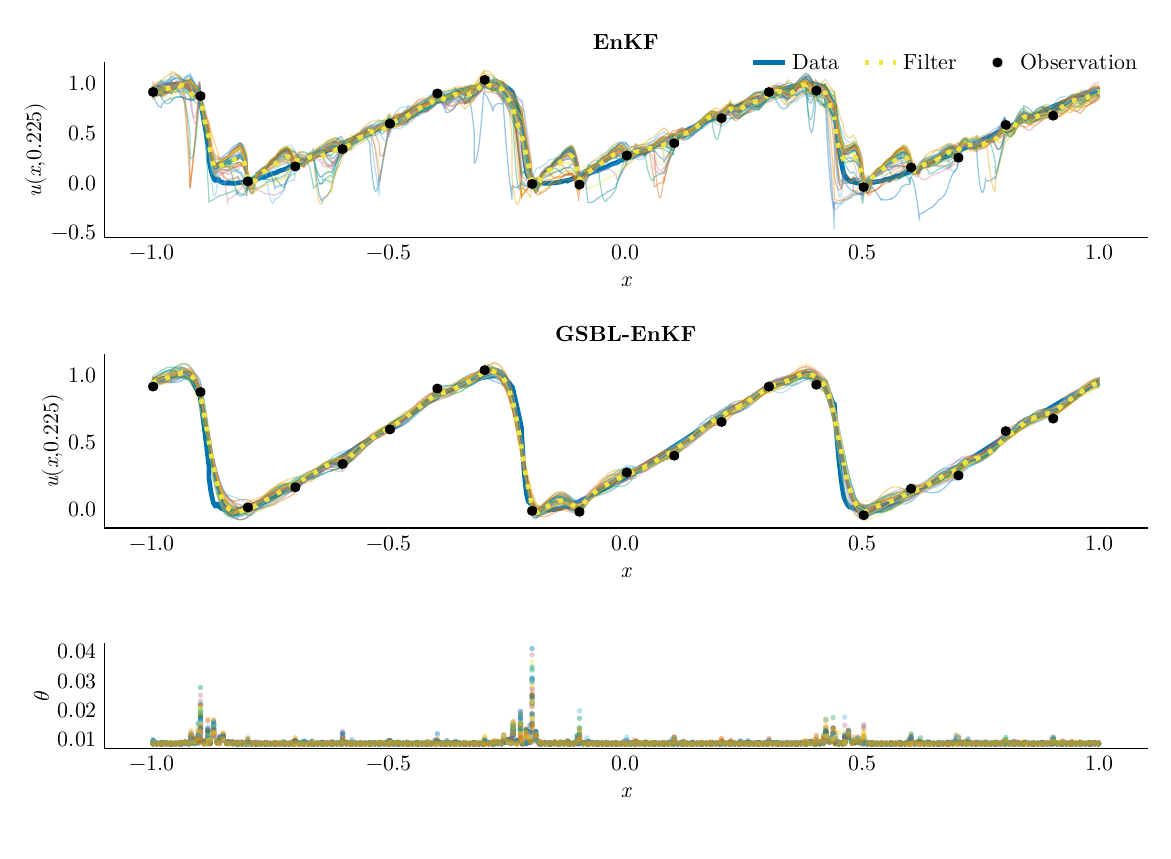}
    \caption{Spatial profiles of the EnKF samples and our proposed GSBL-EnKF samples for the inviscid Burgers equation in \Cref{sub:tests_Burgers} at time $t=0.225$. Different colors represent different ensemble members, and the dashed curve represents the ensemble mean.}
    \label{fig:burgers_profiles}
\end{figure}

\begin{figure}[tb]
    \centering
    \includegraphics[width=0.99\linewidth,clip,trim={0.1cm 0.45cm 0.1cm 0.1cm}]{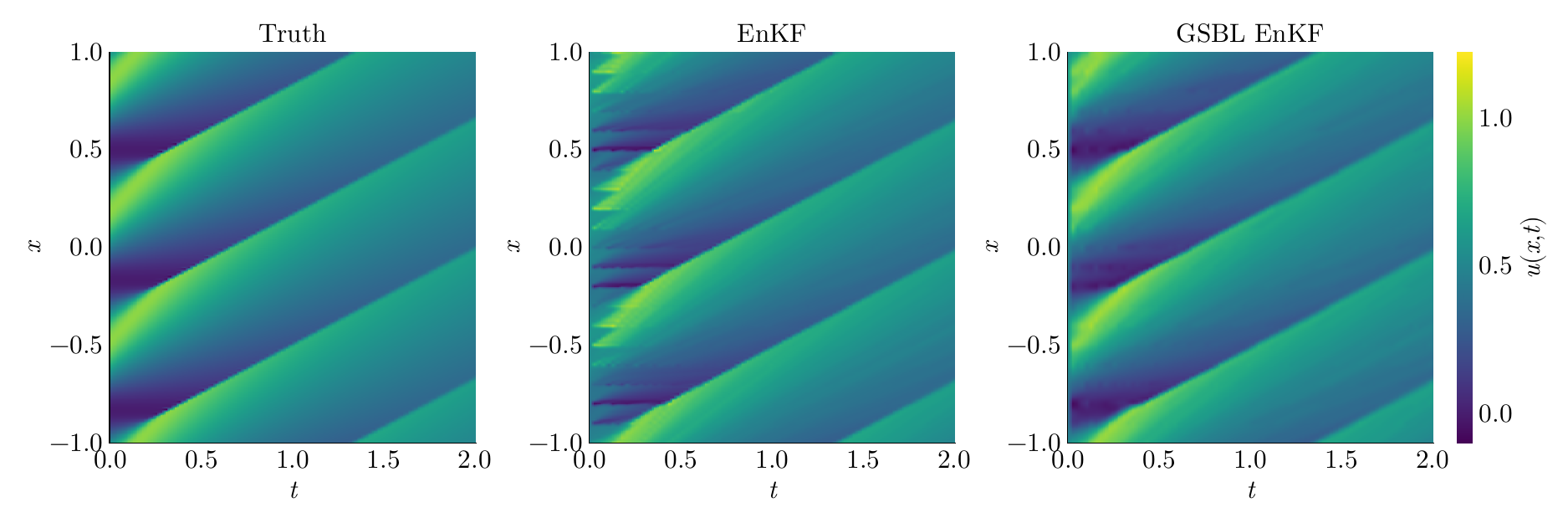}
    \caption{
        {Spatiotemporal mean profiles of the EnKF and our proposed GSBL-EnKF with ensemble size $P=50$ for the inviscid Burgers equation in \Cref{sub:tests_Burgers}.}
    }
    \label{fig:burgers_heatmap}
\end{figure}

We again show the profile of the ensembles of the EnKF and our proposed GSBL-EnKF for a given time in \Cref{fig:burgers_profiles} and a space-time plot of the average ensemble members in \Cref{fig:burgers_heatmap}.
We observe that, while both ensembles tend to improve toward the end of the simulation, the GSBL-EnKF method tends to ``forget'' the erroneous initialization earlier.
Moreover, the GSBL-EnKF suffers less from numerical artifacts from the observation locations, which are located precisely at the horizontal discontinuities prominent at early times in the EnKF panel.

\begin{figure}[tb]
    \centering
    \includegraphics[width=0.8\linewidth,clip,trim={0.1cm 0.15cm 0.1cm 0.1cm}]{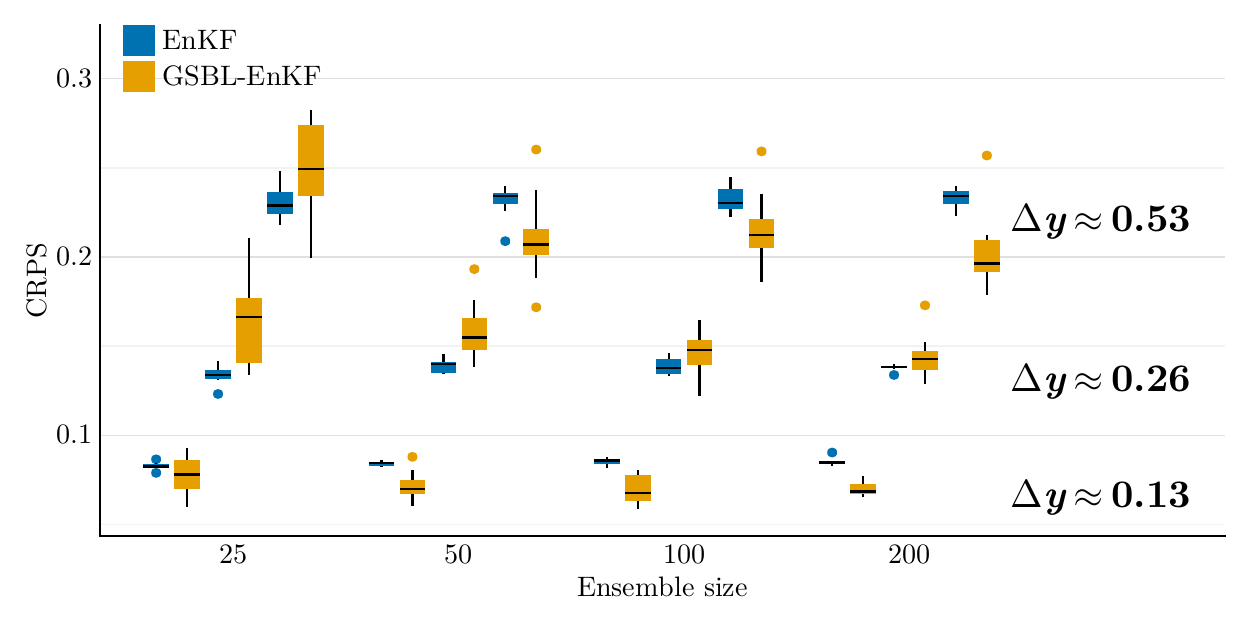}
    \includegraphics[width=0.8\linewidth,clip,trim={0.1cm 0.15cm 0.1cm 0.1cm}]{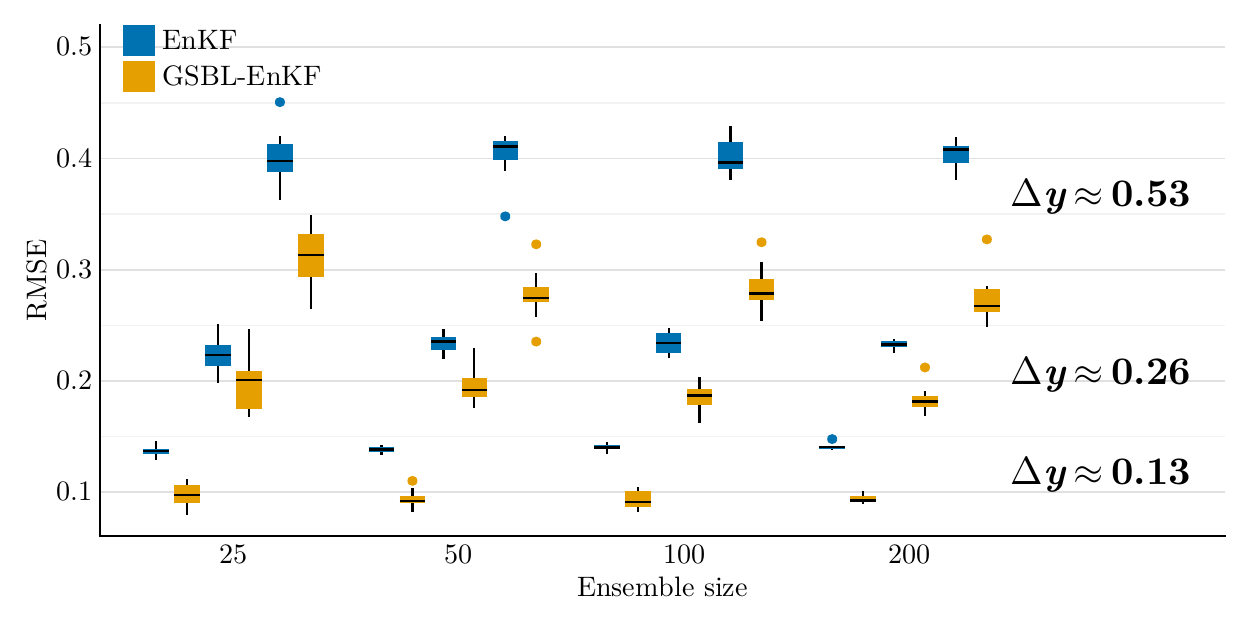}
    \caption{
        The CRPS (top) and RMSE (bottom) of the EnKF and our proposed GSBL-EnKF for different grid spacings and ensemble sizes across independent trials of the inviscid Burgers equation in \Cref{sub:tests_Burgers}.
    }
    \label{fig:burgers_convergence}
\end{figure}

We quantify the improvement of our proposed GSBL-EnKF over the EnKF in \Cref{tab:burgers_improvement_proportion} for different ensemble sizes.
Each count is taken over thirty independent trials, where the filters are initialized identically and assimilate identical data for each trial.
We report the number of trials that our GSBL-EnKF variant is better in (i) time-averaged CRPS and (ii) time-averaged RMSE, taken over ten trials and aggregated across three observation spacings.
We see two contrasting qualities from these frequencies:
For a given initialization, our GSBL-EnKF approach tends to improve on the EnKF, especially as the ensemble size grows.
On the other hand, we also observe that our GSBL-EnKF is more prone to ensemble collapse, as indicated by a relatively smaller improvement in CRPS compared to RMSE.
While our GSBL-EnKF improves time-averaged RMSE across virtually all experiments, the method is penalized for being all-too-certain in regions where it is incorrect.
\Cref{fig:burgers_convergence} illustrates the aggregated results.
While the GSBL-EnKF consistently outperforms the EnKF in RMSE, the CRPS reflects the reduced spread of our approach:
For the most dense observation grid ($\Delta x_{\mathrm{obs}} \approx 0.13$), the GSBL-EnKF has a clearly smaller CRPS for all ensemble sizes; for $\Delta x_{\mathrm{obs}} \approx 0.26$, the two methods are comparable; and for the sparsest grid ($\Delta x_{\mathrm{obs}} \approx 0.53$), the GSBL-EnKF has a larger CRPS for the small ensembles ($P \le 50$) and becomes comparable to the EnKF only for $P \ge 100$. That is, the reduced spread is penalized precisely when the observations are too sparse to quickly correct the erroneous initialization.

\begin{table}[ht!]
    \centering
    \caption{
        Across ten realizations and three observation grids $\Delta y \approx (0.134, 0.268, 0.534)$, i.e., thirty independent trials, frequency of improvement from EnKF to our proposed GSBL-EnKF for the inviscid Burgers equation in \Cref{sub:tests_Burgers}.
        A higher value indicates that our GSBL-EnKF outperforms the EnKF.
    }
    \begin{tabular}{@{}r r r@{}}\toprule
        Ensemble size & CRPS & RMSE \\\midrule
         25 & 26\% (08/30) & 90\% (27/30) \\
         50 & 56\% (17/30) & 100\% (30/30) \\
         100 & 70\% (21/30) & 100\% (30/30) \\
         200 & 76\% (23/30) & 100\% (30/30)\\
         \bottomrule
    \end{tabular}
    \label{tab:burgers_improvement_proportion}
\end{table}

\subsection{Euler's equations: Sod shock-tube}
\label{sub:tests_Euler}

We continue our study with the one-dimensional compressible Euler equations
\begin{equation}\label{eq:Euler}
    \partial_t \begin{pmatrix} \rho \\ \rho v \\ \rho e \end{pmatrix}
    + \partial_x \begin{pmatrix} \rho v \\ \rho v^2 + p \\ (\rho e + p) v \end{pmatrix}
    = 0,
\end{equation}
simulating an ideal gas with density $\rho$, velocity $v$, total energy density $\rho e$, and pressure $p = (\gamma - 1) \left( \rho e - \rho v^2 / 2 \right)$, where the ratio of specific heat is chosen as $\gamma = 1.4$.

\begin{figure}[tb]
    \centering
    \includegraphics[width=\linewidth,clip,trim={0.1cm 0.5cm 0.1cm 0.1cm}]{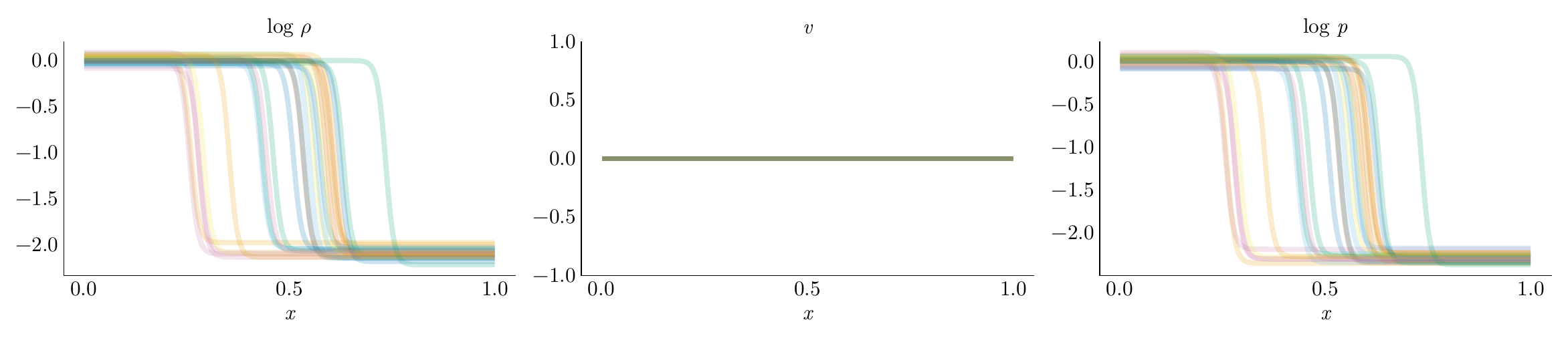}
    \caption{
        {Initialized ensemble for the Sod shock-tube problem for the compressible Euler equations in \Cref{sub:tests_Euler}.}
    }
    \label{fig:euler_init}
\end{figure}

We consider the Sod shock-tube problem~\cite{sodSurveySeveralFinite1978,leveque1992numerical} for the system~\cref{eq:Euler} on $\Omega = (0,1)$ with initial data $\boldsymbol{u}_0 = (\rho,v,p)$ defined as
\begin{equation}\label{eq:sod_ic}
    \boldsymbol{u}_0(x) = \begin{cases}
    (1, 0, 1) & x < 0.5\\
    (0.125, 0, 0.1) & x \geq 0.5.
    \end{cases}
\end{equation}
The Sod shock-tube problem is a standard benchmark for assessing the robustness of numerical schemes, as its solution exhibits three distinct features that probe the propagation of different regions within the system: a rarefaction wave, a contact discontinuity, and a shock discontinuity.
By comparing numerical results with the corresponding analytical solution, one can evaluate how accurately a given method captures and resolves shocks and contact discontinuities, and how well it reproduces the density profile within the rarefaction wave.
This also makes it a relevant setting for our proposed GSBL-EnKF method, as the example tests its ability to handle different types of discontinuities.
To avoid non-positive densities and pressures, we assimilate the transformed primitive variables $\log\rho,v,\log p$, which are shown on this logarithmic scale in any figure.
The sparsifying transform $S$ and the hyperparameters $\boldsymbol{\theta}^{(p)}$ are applied to each of the three components separately.
\Cref{fig:euler_init} illustrates how the different variables are initialized.
We use choices similar to those in~\cite{zhou2026neural}, randomizing the left and right values and the shock location; for numerical stability, however, we smear the shock.\footnote{Recall one motivation for our method is to ensure that we can recover sharp features from smeared initializations.}
More details can be seen in \Cref{app:details}.
As in~\cite{zhou2026neural}, we only collect (noisy) observations of the pressure, i.e., with no information of density or velocity of the fluid, at the sensor locations to better mimic realistic scenarios where not all components may be observed.
Finally, because we do not assume certainty over the boundary conditions, we do not enforce any boundary data in the discretization.

\begin{figure}[tb]
    \centering
    \includegraphics[width=\linewidth,clip,trim={0.1cm 0.45cm 0.1cm 0.1cm}]{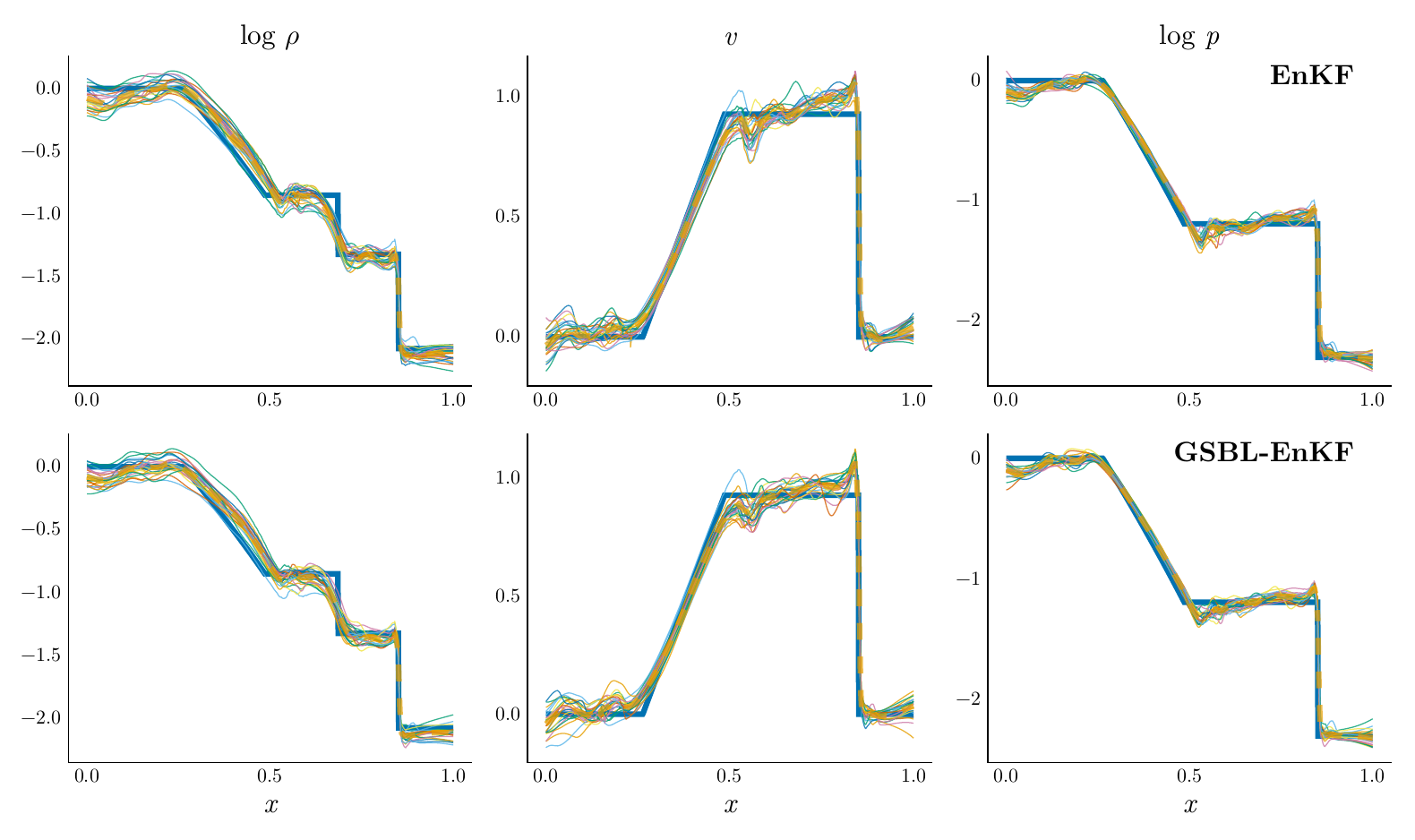}
    \caption{%
        Spatial solution profiles of the EnKF samples and our proposed GSBL-EnKF samples at time $t=0.2$ for the Sod shock-tube problem for the compressible Euler equations in \Cref{sub:tests_Euler}. Different colors represent different ensemble members, and the dashed curve represents the ensemble mean. The solid blue line is the observation-generating simulation, i.e., the data-generating simulation (the reference state) at this time.
    }
    \label{fig:euler_solution}
\end{figure}

\begin{figure}[tb]
    \centering
    \includegraphics[width=0.9\linewidth,clip,trim={0.1cm 0.24cm 0.1cm 0.1cm}]{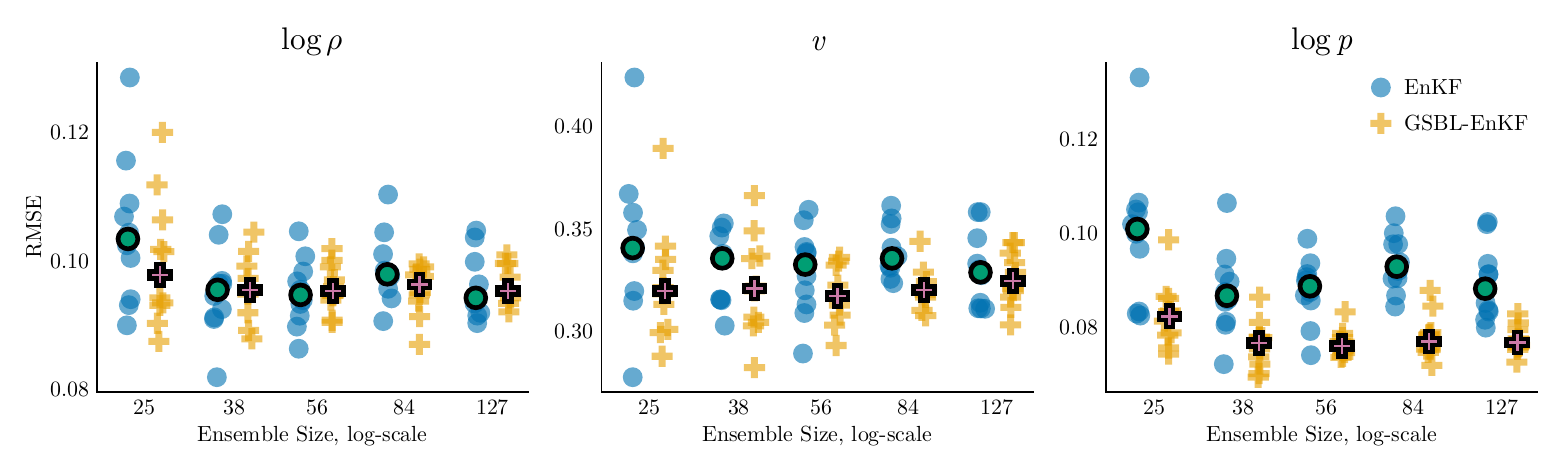}
    \includegraphics[width=0.9\linewidth,clip,trim={0.1cm 0.24cm 0.1cm 0.1cm}]{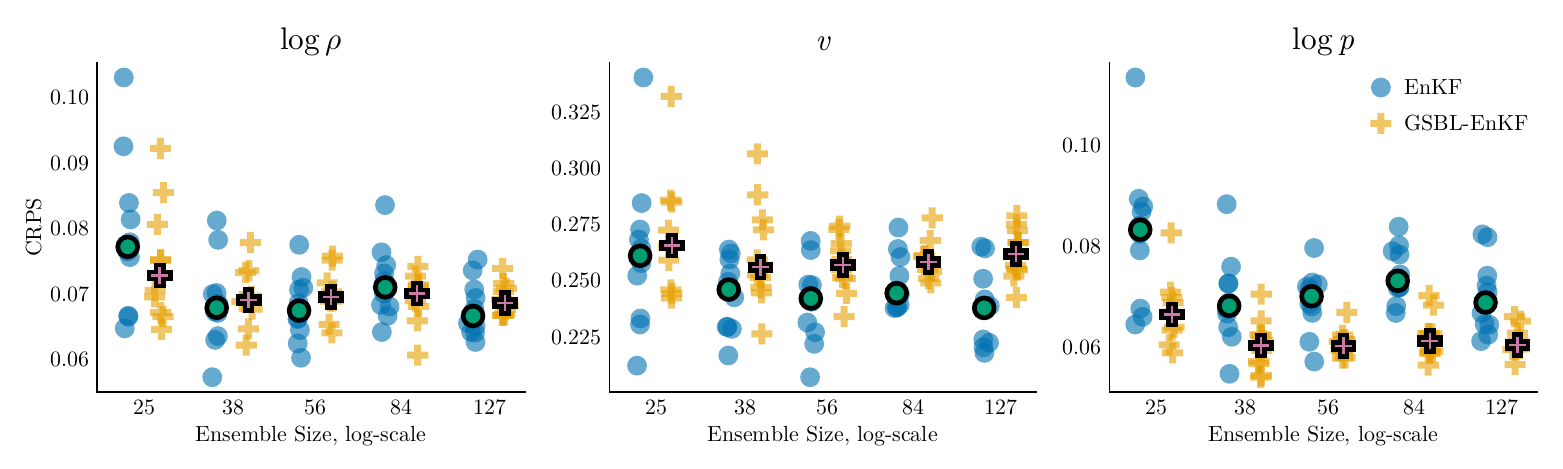}
    \caption{
        {Time-averaged RMSE (top) and CRPS (bottom) for the EnKF and our proposed GSBL-EnKF for different ensemble sizes across ten independent trials of the Sod shock-tube problem for the compressible Euler equations in \Cref{sub:tests_Euler}.
        From left to right, the plots show the time-averaged error for the log-density, velocity, and log-pressure.
        We jitter and offset the abscissae for readability; both methods use the same ensemble sizes.}
    }
    \label{fig:sod_rmse}
\end{figure}

\Cref{fig:euler_solution} reports the spatial solution profiles resulting from our proposed GSBL-EnKF method and the EnKF method at time $t=0.2$.
Notably, although the profiles for both methods appear largely similar in this more challenging test case, we observe that each individual member of our GSBL-EnKF method is smoother in space than in the EnKF method.
Furthermore, \Cref{fig:sod_rmse} shows that the median RMSE of the GSBL-EnKF is slightly smaller than that of the EnKF for all three variables, but the differences are within the spread across trials. 
In terms of CRPS, the GSBL-EnKF improves on the EnKF only for the observed variable, the pressure; for the unobserved density and velocity, its CRPS is comparable or larger, reflecting a reduced ensemble spread in the unobserved components. 
Overall, the improvements are smaller than in the previous test cases. 
We attribute this to the short assimilation window (eight analysis steps) and to the fact that only the pressure is observed.

In fact, while performing our computational experiments, we identified several broader challenges in DA for hyperbolic conservation laws and using high-order methods for their discretization.
In \Cref{sub:problems}, we will lay out a few of the more pressing issues and discuss their implications on what we consider the most relevant direction for future research on DA for hyperbolic conservation laws.

\section{{Concluding remarks}}
\label{sec:summary}

We conclude by summarizing the paper (\Cref{sub:summary}) and discussing open problems together with directions for future research (\Cref{sub:problems}).

\subsection{Summary}
\label{sub:summary}

We introduce a regularity-informed filtering framework for DA in time-dependent PDEs whose solutions exhibit steep gradients and jump discontinuities, with hyperbolic conservation laws as the guiding application.
The core idea is to augment the optimization formulation~\cref{eq:EnKF_analysis2} of the EnKF's analysis step with an edge-preserving spatial regularization term, which we encode using hierarchical GSBL priors.
The resulting GSBL-EnKF method tilts the traditional EnKF prior toward piecewise smooth states:
A conditionally Gaussian prior promotes sparsity of the transformed states $S \mathbf{u}^{(p)}$, while spatially varying hyper-parameters $\boldsymbol{\theta}^{(p)}$, endowed with a generalized gamma hyper-prior, are inferred jointly with the states.
The GSBL-EnKF analysis step corresponds to the MAP estimate of the resulting joint posterior~\cref{eq:posterior}, which we approximate efficiently using a block-coordinate descent algorithm; see \Cref{alg:GSBL_EnKF}.
Each IAS iteration alternates between a regularized Kalman update for the states and explicit, decoupled updates for the hyperparameters.
Only a few iterations are needed in practice ($\nIAS = 2$ in our experiments).
In this sense, the prior encodes our belief about the state's regularity; the optimization formulation in the analysis step ensures that acting on this belief is computationally affordable at each assimilation step.

We demonstrate the GSBL-EnKF on benchmark problems of increasing difficulty: the linear advection equation with piecewise linear data, the inviscid Burgers equation, and the Sod shock-tube problem for the compressible Euler equations.
Compared to the traditional EnKF, the GSBL-EnKF yields ensemble members that better preserve jump discontinuities and are less oscillatory in smooth regions, and it reduces the error of the ensemble mean in nearly all trials.
Qualitatively, it also recovers more quickly from erroneous initializations.
At the same time, the regularized update reduces the ensemble spread. For the more challenging examples, the improvement in CRPS is confined to dense observation grids and larger ensembles, and the GSBL-EnKF is overconfident for sparse observations, small ensembles, and unobserved variables.
For the Sod problem, only the observed pressure is improved consistently in the CRPS.
In addition, the inferred hyperparameters localize the discontinuities of each ensemble member, providing an interpretable indication of possible shock locations and their uncertainty.
While we formalize our approach for the EnKF, hyperbolic conservation laws, and GSBL priors, the underlying principle---incorporating problem-specific regularity assumptions into the analysis step---carries over to other dynamical systems and filtering techniques.
We hope this work paves the way for regularity-informed DA methods that align more closely with the physical properties of complex dynamical systems.

\subsection{Open problems and future work}
\label{sub:problems}

Our computational experiments surface a few open problems.
We believe these are rooted in modeling choices that are not yet well understood and whose resolution matters for DA of hyperbolic conservation laws beyond the specific method proposed here.
We discuss them below, together with several directions for future research that we consider particularly relevant.

\paragraph*{Sparsifying transforms}

We chose the sparsifying transform $S$ to be a weighted second-derivative operator implemented via the spatial discretization of the PDE itself.\footnote{More specifically, transform $S$ maps the state values at the quadrature nodes to the interpolation nodes, applies the second-derivative operator in this space, and then maps back to the quadrature nodes.}
While other structurally-informed methods are also inspired by derivative information, e.g.,~\cite{li2024structurally}, we are not aware of prior work evaluating this derivative through the numerical discretization, whose action on a high-order discontinuous Galerkin or finite-volume approximation can differ markedly from a finite-difference approximation.
The impact of these choices on the induced prior---and thus on the resulting filter---remains to be understood.
More broadly, the GSBL-EnKF framework allows for any suitable sparsifying transform, and the choice should reflect the regularity of the underlying solution.
This is particularly challenging for states with spatially varying smoothness.

\paragraph*{Alternative priors}

We encoded piecewise smoothness using GSBL priors~\cite{calvetti2020sparse,glaubitz2023generalized}, whose conditionally Gaussian structure with generalized gamma hyper-priors enables efficient inference for linear data models with Gaussian noise~\cite{calvetti2020sparsity,xiao2023sequential,lindbloom2024generalized}.
The same structural belief can be promoted by other edge-preserving priors, including (hyper-)Laplace~\cite{figueiredo2007majorization,babacan2009bayesian,krishnan2009fast}, Cauchy~\cite{markkanen2019cauchy,suuronen2022cauchy}, horseshoe~\cite{carvalho2009handling,uribe2023horseshoe,dong2023inducing}, and Besov~\cite{lassas2009discretization,dashti2012besov,lan2023spatiotemporal} priors, to which our general framework extends. A systematic comparison within ensemble filtering would be valuable.

\paragraph*{Nonlinearity introduced by shock capturing}

Semidiscretizations of conservation laws typically rely on shock-capturing schemes when evaluating the forward operator $\Psi$, so that even dynamics that are linear at the PDE level---such as the advection equation in \Cref{sub:tests_linear}---yield a nonlinear forward operator.
Quantifying how this solver-induced nonlinearity affects the theory and expected convergence behavior of the EnKF is an open problem.

\paragraph*{Ensemble spread and variance collapse}

The GSBL-EnKF exhibits a smaller ensemble variance than the traditional EnKF.
For instance, in \Cref{sub:tests_Burgers}, our GSBL-EnKF improved the time-averaged RMSE in virtually all trials yet was penalized in CRPS for being overconfident in regions where it was incorrect.
Increasing the stochastic perturbations (or inflation) could counteract the decreasing spread.
Care must be taken, however, to avoid violating physical constraints, such as the positivity of density in the Euler equations.
At the same time, part of the contraction is physically meaningful:
The dynamics considered here are neither chaotic nor unstable, and shock formation dissipates the mathematical entropy (e.g., the $L^2$ energy for the Burgers equation) and progressively eliminates small-scale information.
For instance, in the periodic Burgers equation, a smooth initial condition evolves toward a constant steady state.
A filter should reflect this problem-intrinsic loss of information.
Disentangling such physically consistent contraction from spurious overconfidence and designing spread-control mechanisms that respect physical constraints remain open problems.

\paragraph*{Mesh-adaptive filtering}

Since many schemes for conservation laws refine the discretization around shocks, future work should consider mesh-sensitive filtering, in which ensemble members may not be well-represented on a shared mesh.
Aligning meshes and shocks across the ensemble points toward filters in which each member carries its own, possibly shock-fitted, discretization.

\paragraph*{Benchmark design}

For deterministic dynamics, such as the ones considered here, the initial condition is often the only source of stochasticity in the simulation.
Unlike chaotic benchmarks such as the Lorenz-63 and -96 systems, which challenge a filter through their sensitivity to inputs, the problems considered here behave benignly away from discontinuities.
Initializing the ensemble via small random perturbations of the true initial condition therefore reveals little about a method's improvement over the EnKF, and the sensitivity to initialization must instead be tested thoroughly.
Our setups thus use deliberately adversarial initializations, as described in \Cref{app:details}. Notably, in the regime of noiseless dynamics, the filtering problem is equivalent to Bayesian inference of the initial condition.

\paragraph*{State representations for inference}

Many papers on Bayesian filtering discuss how to represent a state.
That is, they ask which variables we should use in the EnKF for a given application.
In particular, our examples represent the state at the quadrature points of a discretization.
By contrast, using the interpolation points from our discretization would yield competing spatial values at element interfaces, making localization more difficult.
Returning to recent literature on similar topics,~\cite{zhou2026neural}~uses neural networks to represent the state.
Other latent-space approaches include~\cite{tongLatentAutoencoderEnsemble2026}, which introduces an approach based on autoencoders, or the previously-mentioned~\cite{hansen2024normal}, which performs marginal transformations before assimilating (similar to our approach in \Cref{sub:tests_Euler}).
A more traditional numerical method includes~\cite{aydogduDataAssimilationUsing2019}, which represents uncertainty for moving-mesh problems by projecting to and from a uniform grid.
Finally, a paper that more thoroughly discusses the state representation over classical function approximation classes is~\cite{butlerReparameterizationStatisticalState2012}, which shows examples using polynomial chaos expansions and briefly discusses appropriate approximation classes for functions with discontinuities.
These references are by no means exhaustive---other approaches considering, e.g., projection and spectral collocation, abound---rather, we consolidate a brief set of approaches representing that, in systems with information about state regularity, the choice of what one might represent uncertainty over is \textit{crucial}.

\paragraph*{Preserving invariants}

Combining the proposed approach with invariant-preserving DA methods~\cite{simon2010kalman,janjic2014conservation,wu2019adding,provost2024preserving,subrahmanya2024preserving}, which enforce physical constraints such as positivity and the conservation of mass, energy, or entropy, could enhance the numerical solution's physical fidelity and long-term stability.

\section*{Funding}
JG acknowledges support from the Swedish Research Council (VR) Starting Grant \#2025-05370, the Zenith Career Development Grant \#26.07, and the National Academic Infrastructure for Supercomputing in Sweden (NAISS) grants \#2025/22-1599 and \#2024/22-1207.
DS and YM acknowledge support from the US Department of Energy (DOE), Office of Science, Office of Advanced Scientific Computing Research (ASCR), via the FASTMath6 SciDAC Institute via contract number DE-AC52-07NA27344, under a subcontract from Lawrence Livermore National Laboratory, and via the M2dt MMICC center under award number DE-SC0023187. YM and JG also acknowledge support from the Office of Naval Research, ANSRE MURI, under award N00014-20-1-2595.

\section*{Acknowledgements}
We thank Tongtong Li for insightful discussions and constructive feedback.

\section*{Dedication}
We dedicate this manuscript to the memory of our friend and colleague Mathieu Le Provost.
Mathieu initiated this line of research and brought us together to work on it, but passed away before it could be completed. We hope to honor his efforts and ideas by finishing this manuscript on his behalf.

\bibliographystyle{siamplain}
\bibliography{references}

\appendix
\section{Proof of \texorpdfstring{\Cref{lem:Kalman}}{Lemma \ref{lem:Kalman}}}
\label{app:proofs}

Let $\what{C}$ be SPD.
Because all derivation in this section is given a single state ensemble member $\mathbf{u}$, parameter ensemble member $\boldsymbol{\theta}$, and perturbed observation $\mathbf{b}^{(p)}$, we will abuse notation and drop the superscripts with respect to ensemble index $p$. Now, recall that the GSBL-Kalman update of state $\mathbf{u}$ is equivalent to~\cref{eq:u_update3}, which we write as
\begin{align}\label{eq:proof_Kalman1}
	\left( H^\top  \Gamma^{-1} H + \what{C}^{-1} + S^\top  \Theta^{-1} S \right)&\left( \mathbf{u} - \what{\mathbf{u}} \right) = H^\top  \Gamma^{-1} \left( \mathbf{b} - H \what{\mathbf{u}} \right) - S^\top \Theta^{-1} S \what{\mathbf{u}}.\nonumber
\end{align}
Denoting $H_r = [H; S]$ and $\Gamma_r = \diag( \Gamma, \Theta )$, we express~\cref{eq:proof_Kalman1} more compactly as
\begin{equation}\label{eq:proof_Kalman2}
	\left( H_r^\top \Gamma_r^{-1} H_r + \what{C}^{-1} \right) \left( \mathbf{u} - \what{\mathbf{u}} \right)
	 	= H_r^\top \Gamma_r^{-1} \left(
            \begin{bmatrix}
            		\mathbf{b} \\ \mathbf{0}
            \end{bmatrix}
            - H_r \what{\mathbf{u}}
		\right).
\end{equation}
Note the equivalence
\begin{equation}\label{eq:proof_Kalman25}
	(H_r^\top \Gamma_r^{-1} H_r  + \what{C}^{-1} )^{-1} = \what{C} (H_r^\top \Gamma_r^{-1} H_r \what{C} + I )^{-1}.
\end{equation}
If we left-multiply both sides of~\cref{eq:proof_Kalman2} by $(H_r^\top \Gamma_r^{-1} H_r + \what{C}^{-1})^{-1}$ and use~\cref{eq:proof_Kalman25}, we get that
\begin{equation}\label{eq:proof_Kalman3}
	\mathbf{u} - \what{\mathbf{u}}
	 	= \what{C} \left( H_r^\top \Gamma_r^{-1} H_r \what{C} + I \right)^{-1} H_r^\top \Gamma_r^{-1} \left(
            \begin{bmatrix}
            		\mathbf{b} \\ \mathbf{0}
            \end{bmatrix}
            - H_r \what{\mathbf{u}}
		\right).
\end{equation}
We use the push-through identity $( U V + I )^{-1} U = U ( V U + I )^{-1}$ (see~\cite{henderson1981deriving}). Using $U = H_r^\top$ and $V = \Gamma_r^{-1} H_r \what{C}$ implies the identity:
\begin{equation}
	( H_r^\top \Gamma_r^{-1} H_r \what{C} + I )^{-1} H_r^\top = H_r^\top ( \Gamma_r^{-1} H_r \what{C} H_r^\top + I )^{-1}.
\end{equation}
Thus, the push-through identity shows that~\cref{eq:proof_Kalman3} is equivalent to:
\begin{equation}\label{eq:proof_Kalman4}
	\mathbf{u}
	 	= \what{\mathbf{u}} + \what{C} H_r^\top \left( \Gamma_r^{-1} H_r \what{C} H_r^\top + I \right)^{-1} \Gamma_r^{-1} \left(
            \begin{bmatrix}
            		\mathbf{b} \\ \mathbf{0}
            \end{bmatrix}
            - H_r \what{\mathbf{u}}
		\right).
\end{equation}
Finally, recall the definition $K_r = \what{C} H_r^\top \left( \Gamma_r^{-1} H_r \what{C} H_r^\top + I \right)^{-1} \Gamma_r^{-1}$, whose use in~\cref{eq:proof_Kalman4} precisely yields the assertion.

\section{Implementation Details}
\label{app:details}
Here, we describe our many numerical choices.

\subsection{Solver discretizations}

We use a nodal discontinuous Galerkin semidiscretization based on Gaussian summation-by-parts operators, as implemented in~\cite{ranocha2022adaptive}. For each example, we use a second-order spatial discretization with 100 uniformly sized elements. We include other details in \Cref{tab:solver_discretizations}.

\subsection{Error metrics}

We define our error metrics as the time-averaged root-mean-squared error (RMSE) and the continuous ranked probability score (CRPS):
\begin{gather}
    \mathcal{L}_{\bullet}\big(\mathbf{u}^{(1:P)}_{0:J}\big) \coloneq \frac{1}{J+1}\sum_{j=0}^J \mathcal{L}_\bullet\big(\mathbf{u}_j^{(1:P)},\mathbf{b}_j\big),\\
    \mathcal{L}_{\mathrm{RMSE}}^2(\mathbf{u}^{(1:P)},\mathbf{b}) \coloneq \frac{1}{P}\sum_{p=1}^P \|\mathbf{u}^{(p)} - \mathbf{b}\|^2,\\
    \mathcal{L}_{\mathrm{CRPS}}(\mathbf{u}^{(1:P)},\mathbf{b}) \coloneq \mathcal{L}_{\mathrm{RMSE}}(\mathbf{u}^{(1:P)},\mathbf{b}) - \frac{1}{2}V(\mathbf{u}^{(1:P)}),\\
    V^2(\mathbf{u}^{(1:P)}) \coloneq \frac{1}{P(P-1)} \sum_{p=1}^{P} \sum_{p^\prime=1}^P \|\mathbf{u}^{(p)} - \mathbf{u}^{(p^\prime)}\|^2.
\end{gather}

\subsection{Initialization}
For both \Cref{sub:tests_Burgers,sub:tests_linear} we use smooth, periodic initializations via a stochastic random field defined as
\begin{equation}
    \varphi(x,\omega) = \sqrt{2}\,\mathrm{Re}\left[\sum_{k=1}^{N} Z_k(\omega)\exp\left(-0.5k^\alpha + \pi i (k-1) x\right)\right],
\end{equation}
where $Z_k$ are random variables distributed according to a complex standard normal distribution, i.e., $\mathrm{Law}(Z_k) = \mathrm{Law}(A + Bi),$ where $A,B\sim\mathcal{N}(0,1/2).$ Then, in \Cref{sub:tests_linear}, we initialize our ensemble by sampling the random field defined by
\begin{equation}
    u_0^{\mathrm{advection}}(x,\omega) = 0.5 + 0.5\varphi(x,\omega).
\end{equation}
Similarly, in \Cref{sub:tests_Burgers}, we initialize our ensemble by sampling the process defined by
\begin{equation}
    u_0^{\mathrm{Burgers}}(x,\omega) = \varepsilon u_0(x) + (1-\varepsilon)(\varphi(x,\omega) / 3 + 0.5),
\end{equation}
where $u_0(x) = 0.5 + 0.5\sin(3\pi x)$ is the initial condition of the data-generating simulation and $\varepsilon=0.2.$

Finally, for the Sod shock-tube, we initialize via
\begin{gather}
    u_0^{\mathrm{Sod}}(x,\omega) = (u_R(\omega) - u_L(\omega))\ \mathrm{sigmoid}(x - x^*(\omega)),\quad \mathrm{sigmoid}(t) \coloneqq 1/(1 + \exp(\zeta\,t))\\
    u_R\sim\mathcal{N}((0.125, 0, 0.1),\ \mathrm{diag}(0.05, 0.0, 0.05)^2)\\
    u_L\sim\mathcal{N}((1, 0, 1),\ \mathrm{diag}(0.05, 0, 0.05)^2)\\
    x^*\sim\mathrm{Truncated}\left[\mathcal{N}(0.5, 0.125^2), (10^{-3}, 1-10^{-3})\right],
\end{gather}
where the shock parameter $x^*$ is distributed as a truncated normal taking values in the range $0.5 \pm (0.5-10^{-3})$. We set the sigmoid shape parameter, $\zeta$, as $10^2$.

In all our examples, we use a localization matrix $L$ for both methods presented.
We localize the covariance via
\begin{equation}
    [L\odot \what{C}]_{ij} = [\what{C}]_{ij}\ \rho(|x_i - x_j| / \ell),
\end{equation}
for some set length-scale $\ell$, where $\what{C}$ is the empirical covariance matrix of the samples and $\rho$ is the Gaspari--Cohn localization function~\cite{gaspariConstructionCorrelationFunctions1999}. For both filtering methods, we use the common technique of inflation: assuming we are at time $j > 0$ and given the ensemble mean $\boldsymbol\mu \coloneq \frac{1}{P}\sum_{p=1}^{P} \mathbf{u}^{(p)},$ we update each ensemble member before each forecast step according to
\newcommand{\binfl}{\beta_{\mathrm{infl}}}
\begin{equation}
    \mathbf{u}^{(p)}_{\mathrm{infl}} = \mathbf{u}^{(p)} + \binfl(\mathbf{u}^{(p)} - \boldsymbol{\mu}),
\end{equation}
where we set the inflation parameter in our experiments to $\binfl = 0.02$.

For the GSBL method, we use $n_{\mathrm{IAS}} = 2$ optimization iterations.
For the generalized gamma hyperprior, we consistently set the rate as $r = 0.5$ and the shape parameter as $\beta = 5.95$ (in the parametrization of \Cref{subsub:hyper_prior}) across experiments.
That is, $\tau = r\beta - 3/2 = 1.475 > 0$, so that condition (ii) after~\cref{eq:eq:theta_update_IVP} is satisfied and the initial value in~\cref{eq:eq:theta_update_IVP} is $\varphi(0) = (\tau/r)^{1/r} \approx 8.70$.
We remark that the update~\cref{eq:theta_update3} depends on $\beta$ only through the initial value $\varphi(0) = (\tau/r)^{1/r}$ of the IVP~\cref{eq:eq:theta_update_IVP}. For $r = 1/2$ (more generally, whenever $1/r$ is an even integer), the shape parameters $\beta$ and $\beta' = 3/r - \beta$ therefore yield the same $\boldsymbol{\theta}$-update, since $(\tau/r)^{1/r} = (-\tau/r)^{1/r}$; for instance, $\beta = 0.05$ produces the same update as $\beta = 5.95$, although only the latter satisfies condition (ii) and hence corresponds to a well-posed MAP problem.
We also consistently initialize every ensemble hyperparameter vector as $\theta^{(p)}_k = 1$. Finally, we include problem-dependent (hyper-)parameters in \Cref{tab:solver_discretizations}.

\begin{table}[ht!]
    \centering
    \caption{Parameters for \Cref{sec:tests}.}
         \begin{tabular}{@{} r r r r @{}}\toprule
        & \multicolumn{3}{c}{Examples}\\ \cmidrule{2-4}
        \multicolumn{1}{c}{Parameter} & Advection & Burgers & Euler \\\midrule
        \multicolumn{4}{c}{Discretization hyperparameters}\\\midrule
        Observation interval  & 0.5   & 0.025 & 0.025 \\
        Final time              & 20    & 2     & 0.2\\
        \midrule
        \multicolumn{4}{c}{Filtering hyperparameters}\\\midrule
        Initialization $\alpha$ & 0.8   & 0.7   & ---  \\
        Localization $\ell$     & 0.025 & 0.015 & 0.1  \\
        State noise $\sigma_x$  & 0     & 0.05  & 0    \\
        \midrule\multicolumn{4}{c}{GSBL hyperparameters}\\\midrule
        Reference $\vartheta$   & $10^{-1}$ & $10^{-3}$ & $10^{-3}$\\
        EnKF penalty $\lambda$ & 20 & 5 & 1\\
        \bottomrule
    \end{tabular}
    \label{tab:solver_discretizations}
\end{table}

\end{document}